\documentclass[12pt,leqno]{amsart}
\usepackage{amsmath,amssymb,amsthm}
\usepackage[pdftex]{graphicx}

\begin{document}

\title[Irreducible representations of $WB_n$]{Irreducible representations  of  welded braid group}

\author{Inna Sysoeva}

\email{inna.sysoeva.math@gmail.com}
\date{December 30, 2024}

\begin{abstract} 

In this paper we study irreducible matrix representations of the welded braid group $WB_n$, also known as the group of conjugating automorphisms of a free group $F_n.$ We prove that $WB_n$ has no irreducible representations of dimension $r,$ where $2\leqslant r\leqslant n-2$ for $n\geqslant 5.$ We give complete classification of all extensions of  irreducible representations of the braid group $B_n$ to the welded braid group $WB_n$ of dimensions $n-1$ (for $n\geqslant 7$) and $n$ (for $n\geqslant 7,$ $n\neq 8$). Classification of all extensions of  the irreducible $n-1-$ dimensional reduced Burau representation is given for $n\geqslant 5,$ $n\neq 6.$ A new one-parameter family of $n-1-$ dimensional irreducible representations of $WB_n$ is discovered. Classification of all extensions of the  irreducible $n-$dimensional standard representation (also known as Tong-Yang-Ma representation) is given for all $n\geqslant 3.$ New families of the 3-dimensional irreducible representations of $WB_3$ are discovered.

\end{abstract}

\maketitle

\section{Introduction} 

Welded braid group $WB_n$ is a generalization of a classical braid group $B_n,$ for which some crossings are allowed to be “welded”. From the group-theoretical point of view, $WB_n$ is finitely presented  by braid-like and permutation-like generators and relations (see Definition 2.1).\\

In their paper \cite{FRR}, Fenn, Rimányi and  Rourke have proven that the welded braid group is isomorphic to a subgroup of automorphisms of a free group, generated by the braid group and the symmetric group. More specifically, let $F_n$  be the  free group  with  generators $x_1, \dots, x_n,$ ($n\geqslant 2$), and $Aut (F_n)$ be the group of its (right) automorphisms. 

Then the subgroup of $Aut(F_n)$ generated by the automorphisms $\sigma_1, \dots \sigma_{n-1},$ given by\\

$\sigma_i:\left\{ \begin{array}{lll}
	x_i &\to&x_{i+1}\\
	
	x_{i+1}& \to& x_{i+1}^{-1}x_{i}x_{i+1}\\
	x_{j}& \to&x_{j}, j\neq i,i+1
\end{array} \right.$ \hskip 2cm for $ \ i=1, \dots n-1,$\\

and by the automorphisms $\alpha_1, \dots, \alpha_{n-1},$  given by \\

$\alpha_i:\left\{ \begin{array}{lll}
	x_i &\to&x_{i+1}\\
	
	x_{i+1}& \to&x_{i}\\
	x_{j}& \to&x_{j}, j\neq i,i+1
\end{array} \right.$ \hskip 2cm for $ \ i=1, \dots n-1.$\\

is isomorphic to $WB_n.$

The  automorphisms $\alpha_1, \dots, \alpha_{n-1}$ generate the subgroup of $WB_n$ isomorphic to the symmetric group $S_n$ on $n$ elements, and the automorphisms $\sigma_1, \dots \sigma_{n-1}$ generate the subgroup of $WB_n$ isomorphic to the classical braid group $B_n$ on $n$ strings (see, for example, E.Artin,  \cite{Artin}). \\

It also should be mentioned that the welded braid group  appears in the literature as {\it the group of conjugating automorphisms of a free group.} An automorphism from $Aut (F_n)$ is called a conjugating automorphism
if it maps every generator $x_i$ into a word of the form \\$W_i^{-1}( x_1 , \dots x_n)x_{\pi(i)} W_i( x_1 , \dots x_n),$ where $W_i( x_1 , \dots x_n)\in F_n$ and \\$\pi\in S_n.$ It has been proven by Savushkina in \cite{Sav}, Theorem 1, that the group of conjugating automorphisms is, in fact, the welded braid group, presented by the generators and relations given in Definition 2.1. \\

The author's interest in the irreducible matrix representations of the welded braid group is  motivated by the following discussion.\\

Consider the chain of subgroups $B_n\leqslant WB_n\leqslant Aut(F_n).$ 

The linearity of the braid group $B_n$ was proven by Bigelow in \cite{B} and, simultaneously and independently, by Krammer in \cite{K}. 

It has been proven by  Dyer, Formanek and Grossman in \cite{DFG} that for $n=2$ the automorphisms group $Aut(F_2)$ of the free group  $F_2$ is linear if and only if the braid group on 4 strings $B_4$ is linear. Together with the results of Bigelow \cite{B} and Krammer \cite{K}, this proves that $Aut(F_2)$ is linear, and, hence, its  subgroup  $WB_2$ is linear as well. 

The situation is different, however, for $n\geqslant 3.$ In their paper \cite{FP} Formanek and Procesi had proven that the automorphisms group $Aut(F_n)$ is not linear for $n\geqslant 3.$ It is still an open question whether the welded braid group $WB_n$ is linear or not for $n\geqslant 3.$ It is not surprising, in view of the following easy observation (Theorem 1.1), that the attempts by some authors to extend the  Lawrence–Krammer–Bigelow representation  of the braid group $B_n$ to $WB_n$ led to non-faithful representation, since, as was proven by Zinno in \cite{Zinno},  it is  irreducible for generic values of the parameters.\\

{\bf Theorem 1.1.} Let $\rho:WB_n\to GL_{r}({\mathbb C})$ be an $r-$dimensional representation of the welded braid  group $WB_n$ for $n\geqslant 3,$ $ r \geqslant 1,$ such that the restriction of $\rho$ onto braid group $\rho|_{B_n}$ is irreducible. Then $\rho$ is {\it not faithful}.

\begin{proof}  By \cite{Chow}, Theorem III, for $n\geqslant 3,$ the center of $B_n,$  $Z(B_n)=\langle \theta ^n\rangle.$ Since $\rho_{|B_n}$ is irreducible, the image of the central element is a constant matrix, $\rho(\theta^n)=\lambda I$ for some $\lambda \in \mathbb{C}^*.$ Hence, $\rho(\theta^n)$ commutes with every $\rho(g), \ g\in WB_n.$ By \cite{Sav}, Theorem 2, $WB_n$ has trivial center. Thus, $\rho $ is not faithful.
	
\end{proof}

Different representations of the welded braid group $WB_n$ were constructed by several authors (see, for example, \cite{BS}), but  no systematic classification has been given.

The main goal of this paper is to classify  all the extensions to $WB_n$ of the {\it irreducible} representations of $B_n$  of dimension up to $n$ for $n$ large enough. 

Complete classification of the irreducible representations of the braid group $B_n$ of dimension up to $n$ was given in a series of papers \cite{Formanek, S,FLSV}, and it was proven by the author in \cite {S1} that $B_n$ has no irreducible representations of dimension $n+1$ for $n\geqslant 10.$ For $n$ large enough,  all irreducible representations of $B_n$ of dimension up to $n$ are either one-dimensional, or have dimensions $n-2,$ $n-1$ or $n.$ 

 We prove that for $n\geqslant 5$ the welded braid group $WB_n$ has no irreducible representations of dimension $r,$ where  $2\leqslant r\leqslant n-2$ (Theorem 6.1), hence no irreducible representation of dimension $n-2$ can be extended to $WB_n.$ 
 
 To classify $n-1-$dimensional extensions of the irreducible representations, first, in Chapter 3, we classify all extensions  of the irreducible specializations of the reduced Burau representation (see Definition 2.3) for $n\geqslant 5, n\neq 6.$ We prove (see Theorem 3.25) that for these values of $n$ every such extension is equivalent to a tensor product of a specific one-dimensional representation $X_n(1,k)$ of $WB_n,$ $k=\pm 1,$ and the specialization of one of the two non-equivalent one-parameter families, $\widetilde{\beta_n }(t)$ and $\widehat{\beta_n }(t)$ (see Definition 3.1). The representation family equivalent to $\widetilde{\beta_n }(t)$ has appeared in \cite{BS}, while the family $\widehat{\beta_n }(t)$ has not appeared before in the literature, to the best of our knowledge. The classification of all $n-1-$dimensional extensions of the irreducible representations for $n\geqslant 7$ is given in Theorem 6.2. We prove that every such extension is equivalent to a tensor product of a one-dimensional representation of $WB_n$ and a  specialization of either $\widetilde{\beta_n }(t)$ or  $\widehat{\beta_n }(t).$ The classification for $n\leqslant 6$ has a number of exceptional cases, and is outside of the scope of this paper.

To classify $n-$dimensional extensions of the irreducible representations, first, in Chapters 4 and 5, we classify all extensions of the irreducible specializations of the  standard representation, which sometimes is also referred to as Tong-Yang-Ma representation (see Definition 2.3). We prove in Chapter 4, Theorem 4.11, that for $n\geqslant 4$ every such representation is equivalent to a tensor product of a specific one-dimensional representation $X_n(1,k)$ of $WB_n,$ $k=\pm 1,$ and the specialization of a two-parameter family of representations $\widetilde{\tau}_n (t,q).$ This family has appeared before in \cite{BS}. In Chapter 5 we consider case $n=3,$  where, in addition to the representations described above, there is a number of exceptions (see Theorem 5.9). These exceptional representations have not appeared before, to the best of our knowledge. We decided to include case $n=3$ because it was possible to extend the methods of Chapter 4 to get the classification in this case, and we wanted to illustrate what kind of exceptions one can expect for the smaller values of $n,$ not included in the general classification theorem. The classification of all $n-$dimensional extensions of the irreducible representations for $n\geqslant 7, n\neq 8 $ is given in Theorem 6.3. We prove that every such extension is equivalent to a tensor product of a one-dimensional representation of $WB_n$ and a specialization of $\widetilde{\tau}_n (t,q).$

We would like to point out that even though we give in this paper the complete classification of the extensions of the {\it irreducible} representations of dimension up to $n$ for $n$ large enough, the full classification of the irreducible representations of $WB_n$ of dimension greater or equal to $n-1$ is far from complete, as the Example 1.2 shows.  This classification problem  is outside of the scope of this paper and remains open.\\

{\bf Example 1.2.} We give an example of an {\it irreducible} representation of $WB_n,$ such that  restrictions of this representation onto both $B_n$ and $S_n$ are {\it reducible.} 

For any $n\geqslant 3,$ for any $\lambda \in \mathbb{C}^*,$ such that $\lambda\neq 1,$ consider  $\widetilde{\tau}_n (1,\lambda)$ (see Definition 4.1). It is easy to see that this representation of $WB_n$ is irreducible, while its restrictions on $B_n$ and $S_n$ are  reducible.

 More specifically, on the braid-like generators $\sigma_i, i=1,2,\dots ,n-1,$ this representation is defined by \\$\widetilde{\tau}_n (1,\lambda)(\sigma_{i})=I_{i-1}\oplus \left( \begin{array}{cc}0&1\\1&0\end{array}\right) \oplus I_{n-i-1},$ and its restriction onto $B_n$ is equal to the $n-$dimensional natural permutation representation of the symmetric group $S_n.$ So, it is reducible, and has two invariant subspaces: the one-dimensional invariant subspace $U_1,$  spanned by the vector $(1,1,\dots ,1)^T,$ and the $n-1-$dimensional invariant subspace $U_2,$ consisting of those vectors  whose coordinates add up to zero.

On the permutation-like generators $\alpha_i, i=1,2,\dots ,n-1,$ the representation $\widetilde{\tau}_n (1,\lambda)$ is defined by \\$\widetilde{\tau}_n (1,\lambda)(\alpha_{i})=I_{i-1}\oplus \left( \begin{array}{cc}0&1/\lambda\\\lambda&0\end{array}\right) \oplus I_{n-i-1}.$ The restriction onto $S_n$ is an  $n-$dimensional representation, which is equivalent to the same natural permutation representation, so  is reducible. 

However, neither $U_1$ nor $U_2$ are invariant under $\widetilde{\tau}_n (1,\lambda)(\alpha_{i}),$ hence $\widetilde{\tau}_n (1,\lambda)$ is irreducible.\\

{\bf Conventions and assumptions:} \\
We will assume throughout the paper that $n\geqslant 3,$ unless specified otherwise.\\
The group multiplication is read {\it left to right}.\\
The vectors of the linear $r-$dimensional space are written as {\it column } vectors.

\section{Notations and Preliminary Results}

{\bf Definition 2.1.} Welded braid group on $n$ strings, $WB_n$ is an abstract group, generated by $\alpha_1, \dots, \alpha_{n-1},\sigma_1, \dots \sigma_{n-1}$ with defining relations:

$\alpha_i^2=1,$ $1\leqslant i\leqslant n-1$\hfill (1)

$\alpha_i \alpha_j=\alpha_j\alpha_i,$ $|i-j|\geqslant 2,$ \hfill(2)

$\alpha_i \alpha_{i+1}\alpha_i=\alpha_{i+1}\alpha_i \alpha_{i+1},$ $1\leqslant i\leqslant n-2$ \hfill(3)

$\sigma_i \sigma_j=\sigma_j\sigma_i,$ $|i-j|\geqslant 2,$ \hfill(4)

$\sigma_i \sigma_{i+1}\sigma_i=\sigma_{i+1}\sigma_i \sigma_{i+1},$ $1\leqslant i\leqslant n-2$ \hfill(5)

$\alpha_i \sigma_j=\sigma_j\alpha_i,$ $|i-j|\geqslant 2,$ \hfill(6)

 $\alpha_i\alpha_{i+1}\sigma_i = \sigma_{i+1}\alpha_i  \alpha_{i+1},$ $1\leqslant i\leqslant n-2$\hfill(7)

$\sigma_i \sigma_{i+1}\alpha_i=\alpha_{i+1}\sigma_i \sigma_{i+1},$ $1\leqslant i\leqslant n-2$\hfill(8)\\

The generators $\alpha_1, \dots, \alpha_{n-1}$ together with permutation group relations (1)--(3) generate a subgroup of $WB_n,$ isomorphic to the symmetric group $S_n;$ the generators $\sigma_1, \dots \sigma_{n-1}$ together with braid relations (4)--(5) generate a  subgroup of $WB_n,$ isomorphic to the braid $B_n.$ We will refer to the relations (6)--(8) as mixed relations.\\

{\bf Lemma 2.2.} Let $WB_n$   be the welded braid group,  $n\geqslant 3.$ Let $\theta=\sigma_1  \sigma_2 \dots  \sigma_{n-1}.$ Then:\\
(a)  $\sigma_i \alpha_{i+1}\alpha_i=\alpha_{i+1}\alpha_i \sigma_{i+1},$ $1\leqslant i\leqslant n-2;$\\
(b) $\sigma_{i+1}=\theta\sigma_{i}\theta^{-1}$ for $1\leqslant i\leqslant n-2;$\\
(c) $\alpha_{i+1}=\theta\alpha_{i}\theta^{-1}$ for $1\leqslant i\leqslant n-2;$\\
(d) $(\alpha_i\alpha_{i+1})^3=1$ for $1\leqslant i\leqslant n-2.$

\begin{proof} (a) immediately  follows from group identities (1) and (7);\\
(b)  is a well-known identity in  $B_n$  (see, for example, \cite{Chow}, Equation (6));\\
(c) follows from group identities (8) and (6);\\
(d) is a well-known identity in $S_n.$\\

\end{proof}

\vskip 0.5cm

Let $t$ be an indeterminate, and let $\mathbb{C}[t ^{\pm 1}]$ be a Laurent polynomial ring over the
complex numbers. We define the following representations of $B_n$ by matrices over
$\mathbb{C}[t ^{\pm 1}].$ (Here $I_{k}$ denotes the $k\times k$ identity matrix.)
\vskip 0.5cm

{\bf Definition 2.3.} One-dimensional representation, or {\it the character,} (see \cite{Formanek}, p. 280) \\ $\chi_n(t):B_n\to \mathbb{C}[t ^{\pm 1}]$ is given by $\chi_n(t)(\sigma_i)=t$ for all $ \ i=1, \dots , n-1.$  \\

The {\it reduced Burau} representation  (see \cite{Formanek}, p. 281)\\ $\beta_n (t) : B_n\to GL_{n-1}(\mathbb{C}[t ^{\pm 1}])$ is defined by
\vskip 0.5cm
\noindent
	$\beta_n (t)(\sigma_1)={\small \left( \begin{array}{rcc|c}
	-t&&0&\\
	&&&0\\
	1&&1&\\
	\hline
	&&&\\
	&0&&I_{n-3}\\
	&&&
\end{array} \right)},  \beta_n (t)(\sigma_{n-1})={\small \left( \begin{array}{c|rrr}
	&&&\\
	I_{n-3}&&0&\\
	&&&\\
	\hline
	&&&\\
	&1&&t\\
	0&&&\\
	&0&&-t
\end{array} \right)},$   \\  

\vskip .5cm
\noindent
$\beta_n (t)(\sigma_i)={\small \left( \begin{array}
	{c|rrr|c}
	&&&&\\
	I_{i-2}&&0&&0\\
	&&&&\\
	\hline 
	
	&1&t&0&\\
	0&0&-t&0&0\\
	&0&1&1&\\
	\hline
	&&&&\\
	0&&0&&I_{n-i-2}\\
	&&&&
\end{array} \right)} ,\ \ \ i=2,\dots n-2.$ 

\vskip .5cm

The {\it  standard representation} (see \cite{S}, Definition 6)\\
$\tau_n(t):B_n \to GL_n ({\mathbb{C}}[t ^{\pm 1}])$ is defined by

\vskip 0.5cm
$\tau_n(t) (\sigma_i)={\small \left( \begin{array}{c|ccc|c}
	&&&&\\
	I_{i-1}&&0&&0\\
	&&&&\\
	\hline 
	
	&0&&t&\\
	0&&&&0\\
	&1&&0&\\
	\hline
	&&&&\\
	0&&0&&I_{n-1-i}\\
	&&&&
\end{array} \right)},$
\vskip .5cm
\noindent
for all $i=1,2,\dots, n-1.$ 

\vskip 0.5cm

Specializing $t$ to a non-zero complex number $z$ gives complex one-dimensional, $(n-1)-$dimensional and $n-$dimensional representations $\chi_n(z),$ $\beta_n (z)$ and $\tau_n(z)$ respectively.

\vskip 0.5cm

{\bf  Remark  2.4.} The representation $\beta_n(t)$ is written here in a slightly  different basis than in \cite{Formanek} for the convenience of the calculations, as will be seen in Chapter 3. It is obvious, however, that our definition is equivalent to the one in \cite{Formanek}; it can be done, for example, by replacing all even-numbered vectors in the basis by their opposites.

\vskip 0.5cm

{\bf Theorem 2.5.} 1) (\cite{Formanek}, Lemma 6) For $z \in  {\mathbb{C}}^*$, the  
representation $\beta_n (z): B_n\to GL_{n-1}({\mathbb{C}})$ is irreducible  if and only if $z$ is not a root of $P_n(t)=1+t+t^2+\dots +t^{n-1}.$\\
2) (\cite{S}, Lemmas 5.3 and 5.4)  For $z \in  {\mathbb{C}}^*$, the 
representation\\ $\tau_n (z): B_n\to GL_{n}({\mathbb{C}})$ is irreducible  if and only if $z\neq 1.$

 \vskip 0.5cm
 
 In Chapter 3, where we will study the extensions of the \\$(n-1)-$dimensional {\it irreducible}  reduced Burau representation, it will be convenient to consider two different bases of the $(n-1)-$dimensional space $V.$
 
 Whenever we would like to emphasize that the matrix of an operator $R$ is written in a specific basis $\mathcal{B},$ we will denote it by $(R)_{[\mathcal{B}]}.$\\
 
 More specifically, let $z \in  {\mathbb{C}}^*,$ such that $P_n(z)\neq 0,$ so that the specialization of the reduced Burau representation $\beta_n(z)$ is irreducible. Let $\mathcal{W}=\{w_1, \dots ,w_{n-1}\}$ be a basis of $V,$ such that in this basis the matrices of $\beta_n(z)$ are given by\\
 
 	$(\beta_n (z)(\sigma_1))_{[\mathcal{W}]}=\left( \begin{array}{rcc|c}
 	-z&&0&\\
 	&&&0\\
 	1&&1&\\
 	\hline
 	&&&\\
 	&0&&I_{n-3}\\
 	&&&
 \end{array} \right),$
 
  \vskip .3cm
 
 $ (\beta_n (z)(\sigma_{n-1}))_{[\mathcal{W}]}=\left( \begin{array}{c|rrr}
 	&&&\\
 	I_{n-3}&&0&\\
 	&&&\\
 	\hline
 	&&&\\
 	&1&&z\\
 	0&&&\\
 	&0&&-z
 \end{array} \right),$   \\  
 
 \vskip .3cm
 
 $(\beta_n (z)(\sigma_i))_{[\mathcal{W}]}=\left( \begin{array}
 	{c|rrr|c}
 	&&&&\\
 	I_{i-2}&&0&&0\\
 	&&&&\\
 	\hline 
 	
 	&1&z&0&\\
 	0&0&-z&0&0\\
 	&0&1&1&\\
 	\hline
 	&&&&\\
 	0&&0&&I_{n-i-2}\\
 	&&&&
 \end{array} \right) ,\ \ \ i=2,\dots n-2.$ \\
 
 As one can see from the explicit expressions of these matrices, the kernels of the operators $\beta_n (z)(\sigma_i)-I$ are conveniently described in basis $\mathcal{W}$ (see Lemma 3.4(c)).
 
 We will introduce another basis $\mathcal{V}=\{v_1, \dots ,v_{n-1}\}$ of $V,$ such that in this basis the images of  the operators $\beta_n (z)(\sigma_i)-I$ are conveniently described. Moreover, each of the vectors $v_i$  spans  the one-dimensional space $\beta_n (z)(\sigma_i)-I$ (see Lemma 3.4(a)). In case when $P_n(z)\neq 0,$ these vectors form a basis of $V.$ \\

 Namely,  define vectors $v_i,\ i=1,\dots , n-1$ as follows:\\
 $v_1=-z^{n-2}(z+1)w_1+z^{n-2}w_2,$\\
 $v_i=z^{n-i}w_{i-1}-z^{n-i-1}(z+1)w_i+z^{n-i-1}w_{i+1},\ i=2,   \dots , n-2,$\\
 $v_{n-1}=zw_{n-2}-(z+1)w_{n-1}.$\\

 {\bf Lemma 2.6.} Let  $z\in \mathbb{C^*},$ such that $P_n(z)\neq 0.$ Then the vectors $v_1, \dots ,v_{n-1}$ form a basis of $V,$ and in this basis the matrices of $\beta_n(z)$ are given by\\
 
 $(\beta_n(z)(\sigma_1))_{[\mathcal{V}]}=\left( \begin{array}{rcc|c}
 	-z&&1&\\
 	&&&0\\
 	0&&1&\\
 	\hline
 	&&&\\
 	&0&&I_{n-3}\\
 	&&&
 \end{array} \right),$
 \vskip .3cm
 
 $(\beta_n(z)(\sigma_{n-1}))_{[\mathcal{V}]}=\left( \begin{array}{c|rrr}
 	&&&\\
 	I_{n-3}&&0&\\
 	&&&\\
 	\hline
 	&&&\\
 	&1&&0\\
 	0&&&\\
 	&z&&-z
 \end{array} \right),$   
 
\vskip .3cm
 
 $(\beta_n(z)(\sigma_i))_{[\mathcal{V}]}=\left( \begin{array}
 	{c|rrr|c}
 	&&&&\\
 	I_{i-2}&&0&&0\\
 	&&&&\\
 	\hline 
 	
 	&1&0&0&\\
 	0&z&-z&1&0\\
 	&0&0&1&\\
 	\hline
 	&&&&\\
 	0&&0&&I_{n-i-2}\\
 	&&&&
 \end{array} \right) ,\ \ \ i=2,\dots n-2.$ \\

 \begin{proof}
 	The corresponding change-of-basis matrix is
 	\vskip 0.5cm
 	
 	\noindent
 	$Q(z)={\tiny \left( \begin{array}
 		{ccccccc}
 		
 		-z^{n-2}(z+1)&z^{n-2}&0&&&&\\
 		z^{n-2}&-z^{n-3}(z+1)&z^{n-3}&&&\vdots&\\
 		
 		0&z^{n-3}&-z^{n-4}(z+1)&&&&\\
 		&&&\ddots&&&\\
 		&\vdots&&&z^2&-z(z+1)&z\\
 		&&&&0&z&-(z+1)
 	\end{array} \right)}.$\\
 	
 		\vskip 0.5cm
 	Its determinant $det(Q(z))=(-1)^{n-1}z^{\left(\frac{(n-1)(n-2)}{2}\right)}P_n(z)\neq 0, $ so the vectors $\{v_1, \dots ,v_{n-1}\}$ form a basis  $\mathcal{V}$ of $V.$ \\
 	
 	It is easy to check that for all $i=1, \dots , n-1,$ \\
 	$Q(z)\cdot(\beta_n(z)(\sigma_i))_{[\mathcal{V}]}=(\beta_n(z)(\sigma_i))_{[\mathcal{W}]}\cdot Q(z).$\\
 \end{proof}

 Now we will describe all one-dimensional representations of  $WB_n.$
  \vskip 0.5cm
 
 {\bf Definition 2.7.} For every $ y\in \mathbb{C}^*$ and $k \in \{-1,1\}$  define \\$X_n(y,k):WB_n\to\mathbb{C}^*$ by \\$X_n(y,k)(\sigma_i)=y,$  $X_n(y,k)(\alpha_i)=k,$   $1\leqslant i \leqslant n-1.$ 
 
 \vskip 0.5cm
 {\bf Theorem 2.8.}  1) For every $ y\in \mathbb{C}^*$ and $k \in \{-1,1\},$ $X_n(y,k)$ is a one-dimensional representation of $WB_n.$\\ 2) For distinct pais of $(y,k),$ the representations $X_n(y,k)$ are not equivalent.\\3) The representations $X_n(y,k), \ \ y\in \mathbb{C}^*,\ \  k \in \{-1,1\}$  give complete set of one-dimensional representations of $WB_n.$ \\
\begin{proof}
Parts 1) and 2) immediately follow from Definition 2.7 and the identities (1)-(8).\\
3) Suppose $\rho:WB_n\to\mathbb{C}^*$ is a one-dimensional representation of $WB_n.$ Let $\rho(\sigma_i)=y_i$ and $\rho(\alpha_i)=k_i,$ where $y_i,\ k_i \in \mathbb{C}^*,$ $1\leqslant i \leqslant n-1.$ It immediately follows from identities (1),  (3) and (6) that  $y_i=y_{i+1},$ and $k_i=k_{i+1},$ for all $1\leqslant i \leqslant n-2,$ and $k^2_i=1$ for $1\leqslant i \leqslant n-1.$ Thus, $\rho=X_n(y,k)$ for some values  $y\in \mathbb{C}^*$ and $ k \in \{-1,1\}.$

\end{proof}

\vskip 0.5cm

The following easy lemma shows that any two equivalent extensions of an {\it irreducible} representation are equal.
\vskip 0.5cm

{\bf Lemma 2.9.} Suppose  $\rho_1,\rho_2:G\to GL_r(\mathbb{C}),$ are the $r-$dimensional  representations of the  group $G$,  and let $H$ be a subgroup of $G$ such that the restrictions of $\rho_1$ and $\rho_2$ onto  $H$ are equal, $\rho_1|_{H}=\rho_2|_{H}=\hat{\rho}$. Suppose that $\hat{\rho}$ is {\it irreducible}. 

If  $\rho_1$ and $\rho_2$ are equivalent  then $\rho_1=\rho_2.$

\begin{proof} Suppose $P$ is the change of basis matrix between the equivalent representations $\rho_1$ and $\rho_2,$ \\
$\rho_2(g)=P^{-1}\rho_1(g)P$ for all $g\in G.$ \\ In particular, $P$ commutes with $\hat{\rho}(h)$ for all $h\in H.$ Thus, $P$ is a scalar matrix, $P=cI,$ since for an eigenvalue $c$ of $P,$  the non-zero eigenspace $Ker(P-cI)$ is an invariant subspace of an irreducible representation $\hat{\rho}.$ So, $\rho_1(g)=\rho_2(g)$ for all $g\in G.$

\end{proof}

\vskip 0.5cm
The goal of this paper is to find, up to equivalency, all representations $\rho$ of $WB_n$ of degree $n$ such that the restriction of $\rho$ onto $B_n$ is equivalent to the specialization of the standard {\it irreducible} representation. Clearly, all such extensions are irreducible representations of $WB_n.$
\vskip 0.5cm

{\bf Remark 2.10.} Suppose $\rho:WB_n\to GL_r(\mathbb{C})$ is an extension to the group $WB_n$ of an $r-$dimensional  representation of the braid group $B_n,$  \\$\hat{\rho}:B_n\to GL_r(\mathbb{C}).$ Then:

1) Both representations $X_n(1,k)\otimes \rho,$ ( $k\in\{-1,1\}$) are extensions of  the same representation $\hat{\rho}.$ 
 
 2) $X_n(1,1)\otimes \rho=\rho.$
 
 3) If $\hat{\rho}$ is  {\it irreducible}, then, by Lemma 2.9, the extensions  $X_n(1,1)\otimes \rho$ and $X_n(1,-1)\otimes \rho$ are {\it not equivalent}. \\

\vskip 0.5cm

Let $\rho:WB_n\to GL_r(\mathbb{C})$ be an $r-$dimensional  representation of the group $WB_n.$  Throughout the paper, we will use the following notations:\\

\noindent
$G_i=\rho(\sigma_i)=G_i,$ $i=1, \dots, n-1;$\\
 $A_i=\rho(\alpha_i)=A_i,$ $i=1, \dots, n-1;$\\
 $T=\rho(\theta)=\rho(\sigma_1  \sigma_2 \dots  \sigma_{n-1}).$\\
 
 In this paper we are interested in the classification of the extensions to $WB_n$ of the {\it irreducible} representations of $B_n,$ and since for $n$ large enough all irreducible representations of $B_n$ of dimension up to $n,$ after appropriate tensoring with one-dimensional representation, have small corank (see \cite{S} Definition 2, \cite{Formanek} Theorem 22, and \cite{S} Theorem 6.1), it is convenient to consider the operators\\
 
 \noindent
  $C_i=\rho(\sigma_i)-I, \ i=1, \dots, n-1.$\\
  
  We will show in Chapters 3 and 4, that for these extensions, for $n$ large enough, after appropriate tensoring with one-dimensional representation, the rank of $\rho(\alpha_i)-I$ must be the same as the rank of $ \rho(\sigma_i)-I,$ so 
 it is also convenient  to introduce operators \\
 
 \noindent
 $D_i=\rho(\alpha_i)-I=A_i-I,$ for $ i=1, 2, \dots, n-1.$\\

In the following lemma (Lemma 2.11) we will list some properties of these extensions, which we will use throughout the paper. Each of them can be easily obtained from the  group identities (1)--(8).\\

{\bf Lemma 2.11.} Let $\rho:WB_n\to GL_r(\mathbb{C})$ be an $r-$dimensional  representation of  $WB_n,$ $n\geqslant 3.$ Then:\\
(a) $TC_{i}T^{-1}=C_{i+1},$ $i=1, \dots, n-2 ,$ and $T^2C_{n-1}T^{-2}=C_{1};$\\
(b) $A_iC_j=C_jA_i,$ and  $D_iC_j=C_jD_i$ for $1\leqslant i,j\leqslant n-1,$  $|i-j|\geqslant 2;$\\
(c) $A_iA_{i+1}C_i = C_{i+1}A_i  A_{i+1},$ $1\leqslant i\leqslant n-2;$\\
(d)  $C_iA_{i+1}A_{i}= A_{i+1}  A_{i}C_{i+1},$ $1\leqslant i\leqslant n-2;$\\
(e)  $TA_{i}T^{-1}=A_{i+1},$ and $TD_{i}T^{-1}=D_{i+1},$ $i=1, \dots, n-2 ;$\\
(f)   The only possible complex eigenvalues of $A_i,$ $1\leqslant i\leqslant n-1,$ are $1$ and $-1.$ \\

In addition, if the restriction $\rho|_{B_n}$ of $\rho$ onto subgroup $B_n$ is {\it irreducible}, then:\\
(g) $T^n=c I,$ where $c\in \mathbb{C}^{*};$\\
(h)  $T^2A_{n-1}T^{-2}=A_{1},$ and $T^2D_{n-1}T^{-2}=D_{1}.$\\

\begin{proof}
	Parts (a)-(f) immediately follow from the group identities (1)-(8) and Lemma 2.2. \\
	Part (g): By \cite{Chow}, Theorem III, for $n\geqslant 3,$ the center of the braid group   $Z(B_n)=\langle \theta ^n\rangle.$ Since $\rho|_{B_n}$ is irreducible, the image of the central element is a constant matrix.\\
	Part (h) follows from (g) and (e).\\
\end{proof}

 \section{Extensions of the Burau irreducible representation of dimension $n-1.$ }
 
The goal of the this section is to classify all non-equivalent extensions of the $(n-1)-$dimensional {\it irreducible}  reduced Burau representation of the braid group $B_n$ to the welded braid group $WB_n$ for $n$ large enough. The results of this section hold true for all $n\geqslant 7,$ as well as  for $n=5.$  There are additional exceptional extensions for $n=3, 4, 6,$ and the classification for these values of $n$ is left outside of the scope of this paper.\\

We will start this section with the following definition.\\

{\bf Definition 3.1.} Let $t$ be an indeterminate. For every   $n\geqslant 3,$ define \\

$\widetilde{\beta_n }(t) : WB_n\to GL_{n-1}(\mathbb{C}[t ^{\pm 1}])$  by
\vskip 0.5cm

$\widetilde{\beta_n} (t)(\sigma_1)=\left( \begin{array}{rcc|c}
	-t&&1&\\
	&&&0\\
	0&&1&\\
	\hline
	&&&\\
	&0&&I_{n-3}\\
	&&&
\end{array} \right),\ \ \  \widetilde{\beta_n} (t)(\sigma_{n-1})=\left( \begin{array}{c|rrr}
	&&&\\
	I_{n-3}&&0&\\
	&&&\\
	\hline
	&&&\\
	&1&&0\\
	0&&&\\
	&t&&-t
\end{array} \right),$   \\  

\vskip .5cm

$\widetilde{\beta_n} (t)(\sigma_i)=\left( \begin{array}
	{c|rrr|c}
	&&&&\\
	I_{i-2}&&0&&0\\
	&&&&\\
	\hline 
	
	&1&0&0&\\
	0&t&-t&1&0\\
	&0&0&1&\\
	\hline
	&&&&\\
	0&&0&&I_{n-i-2}\\
	&&&&
\end{array} \right) ,\ \ \ i=2,\dots n-2,$ \\

\vskip .5cm

$\widetilde{\beta_n} (t)(\alpha_1)=\left( \begin{array}{rcc|c}
	-1&&1&\\
	&&&0\\
	0&&1&\\
	\hline
	&&&\\
	&0&&I_{n-3}\\
	&&&
\end{array} \right),\ \ \  \widetilde{\beta_n} (t)(\alpha_{n-1})=\left( \begin{array}{c|rrr}
	&&&\\
	I_{n-3}&&0&\\
	&&&\\
	\hline
	&&&\\
	&1&&0\\
	0&&&\\
	&1&&-1
\end{array} \right),$   \\  

$\widetilde{\beta_n} (t)(\alpha_i)=\left( \begin{array}
	{c|rrr|c}
	&&&&\\
	I_{i-2}&&0&&0\\
	&&&&\\
	\hline 
	
	&1&0&0&\\
	0&1&-1&1&0\\
	&0&0&1&\\
	\hline
	&&&&\\
	0&&0&&I_{n-i-2}\\
	&&&&
\end{array} \right) ,\ \ \ i=2,\dots n-2.$ \\

\vskip .5cm

Define $\widehat{\beta_n }(t) : WB_n\to GL_{n-1}(\mathbb{C}[t ^{\pm 1}])$  by
\vskip 0.5cm

$\widehat{\beta_n} (t)(\sigma_1)=\left( \begin{array}{rcc|c}
	-t&&0&\\
	&&&0\\
	1&&1&\\
	\hline
	&&&\\
	&0&&I_{n-3}\\
	&&&
\end{array} \right),\ \ \  \widehat{\beta_n} (t)(\sigma_{n-1})=\left( \begin{array}{c|rrr}
	&&&\\
	I_{n-3}&&0&\\
	&&&\\
	\hline
	&&&\\
	&1&&t\\
	0&&&\\
	&0&&-t
\end{array} \right),$   \\  

\vskip .5cm

$\widehat{\beta_n} (t)(\sigma_i)=\left( \begin{array}
	{c|rrr|c}
	&&&&\\
	I_{i-2}&&0&&0\\
	&&&&\\
	\hline 
	
	&1&t&0&\\
	0&0&-t&0&0\\
	&0&1&1&\\
	\hline
	&&&&\\
	0&&0&&I_{n-i-2}\\
	&&&&
\end{array} \right) ,\ \ \ i=2,\dots n-2,$ \\

\vskip .5cm

$\widehat{\beta_n} (t)(\alpha_1)=\left( \begin{array}{rcc|c}
	-1&&0&\\
	&&&0\\
	\frac{1}{t}&&1&\\
	\hline
	&&&\\
	&0&&I_{n-3}\\
	&&&
\end{array} \right),\ \ \  \widehat{\beta_n} (t)(\alpha_{n-1})=\left( \begin{array}{c|rrr}
	&&&\\
	I_{n-3}&&0&\\
	&&&\\
	\hline
	&&&\\
	&1&&t\\
	0&&&\\
	&0&&-1
\end{array} \right),$   \\  

$\widehat{\beta_n} (t)(\alpha_i)=\left( \begin{array}
	{c|rrr|c}
	&&&&\\
	I_{i-2}&&0&&0\\
	&&&&\\
	\hline 
	
	&1&t&0&\\
	0&0&-1&0&0\\
	&0&\frac{1}{t}&1&\\
	\hline
	&&&&\\
	0&&0&&I_{n-i-2}\\
	&&&&
\end{array} \right) ,\ \ \ i=2,\dots n-2.$ \\

\vskip .5cm

For every $z\in \mathbb{C}^*$ denote by $\widetilde{\beta_n} (z)$ and $\widehat{\beta_n} (z)$ the maps obtained from
$\widetilde{\beta_n} (t)$ and $\widehat{\beta_n} (t)$   respectively by the specialization $t\mapsto z $ .\\ 

{\bf Lemma 3.2.} Suppose $n\geqslant 3.$\\

1) Both $\widetilde{\beta_n} (t)$ and $\widehat{\beta_n} (t)$ are  representations of $WB_n.$ \\

2) For every $z\in \mathbb{C}^*,$ both  $\widetilde{\beta_n} (z)$ and $\widehat{\beta_n} (z)$ are $(n-1)-$dimensional complex representations of $WB_n.$  If  $P_n(z)\neq 0,$ then the restrictions of both $\widetilde{\beta_n} (z)$ and $\widehat{\beta_n} (z)$  onto the braid group $B_n$ are equivalent to an  {\it irreducible} reduced Burau representation   $\beta_n(z).$ \\

3) For every $z_1,z_2\in \mathbb{C}^*,$ such that $P_n(z_1)\neq 0,$  $P_n(z_2)\neq 0,$ and $z_1\neq z_2:$  \\
a) the representations $\widetilde{\beta_n} (z_1)$ and $\widetilde{\beta_n} (z_2)$ are not equivalent; \\
b) the representations $\widehat{\beta_n} (z_1)$ and $\widehat{\beta_n} (z_2)$ are not equivalent; \\
c) the representations $\widetilde{\beta_n} (z_1)$ and $\widehat{\beta_n} (z_2)$ are not equivalent;\\

4) For every $z\in \mathbb{C}^*,$ such that $P_n(z)\neq 0,$ the representations $\widetilde{\beta_n} (z)$ and $\widehat{\beta_n} (z)$ are not equivalent, except when $z=1,$ in which case  $\widetilde{\beta_n} (1)$ is equivalent to $\widehat{\beta_n} (1).$\\

\begin{proof} 1) It is easy to check that all defining relations (1)--(8) of the group $WB_n$  are satisfied for both $\widetilde{\beta_n} (t)$ and $\widehat{\beta_n} (t).$ \\
	
	2) The first statement immediately follows from part 1). By Theorem 2.5, part 1), if $P_n(z)\neq 0,$ then $\beta_n(z)$ is irreducible, and by definition, ${\widetilde{\beta_n} (z)|}_{B_n}=(\beta_n(z))_{[\mathcal{V}]},$ and   ${\widehat{\beta_n} (z)|}_{B_n}=(\beta_n(z))_{[\mathcal{W}]}.$\\

	3) All three statements follow immediately from part 2) and the fact that  $det(\beta_n(z)(\sigma_1))=-z,$ so the representations $\beta_n(z_1)$ and $\beta_n(z_2)$ are not equivalent for $z_1\neq z_2,$ .\\
	
	4) For $z=1,$ the representations $\widetilde{\beta_n} (1)$ and  $\widehat{\beta_n} (1)$  are equivalent, as can be  seen by using the invertible change-of-basis matrix $Q(1).$ Indeed, for every $i=1, \dots , n-1,$ by Lemma 2.6, \\
	$Q(1)(\widetilde{\beta_n}(1)(\sigma_i))= Q(1)(\beta_n(1))_{[\mathcal{V}]}=(\beta_n(1))_{[\mathcal{W}]} Q(1)= (\widehat{\beta_n}(1)(\sigma_i)) Q(1),$ and it is easy to check that   \\ $Q(1)(\widetilde{\beta_n}(1)(\alpha_i))=(\widehat{\beta_n}(1)(\alpha_i)) Q(1).$\\

	Now consider  $z\neq 1,$ such that $P_n(z)\neq 0.$ First, we will write the representation $\widetilde{\beta_n}(z)$ in a different basis, so that in this basis \\${\widetilde{\beta_n} (z)|}_{B_n}={\widehat{\beta_n} (z)|}_{B_n}.$ Namely, by using the change-of-basis matrix $Q(z)$ and writing the representation $\widetilde{\beta_n} (z)$ in basis $\mathcal{W},$ we obtain that for every $i=1, \dots , n-1,$\\
	$(\widetilde{\beta_n}(z)(\sigma_i))_{[\mathcal{W}]}=(\widehat{\beta_n}(z)(\sigma_i))_{[\mathcal{W}]}=(\beta_n(z)(\sigma_i))_{[\mathcal{W}]}$ (see Lemma 2.6).\\
	
	Writing the matrix of $\widetilde{\beta_n}(z)(\alpha_1)$ in the same basis, we obtain\\
	
	$(\widetilde{\beta_n}(z)(\alpha_1))_{[\mathcal{W}]}=I_{n-1}+$\\
	\vskip 0.3cm
	\noindent
	$+\frac{1}{P_{n}(z)}{\tiny \left( \begin{array}{ccccc}
		-(z+1)(P_{n-1}(z)+1)&(z+1)z(z^{n-2}-1)&(z+1)z^{2}(z^{n-3}-1)&\dots&(z+1)z^{n-2}(z-1)\\\
		&&&&\\
	P_{n-1}(z)+1&-z(z^{n-2}-1)&-z^{2}(z^{n-3}-1)&\dots&-z^{n-2}(z-1)\\
		&&&&\\
		0&0&\dots&0&0\\
	
		\vdots&\vdots&\vdots&\vdots&\vdots\\
		
		0&0&\dots&0&0
	\end{array} \right)}.$ \\

\vskip 0.5cm

It is easy to see that the matrix above is indeed the matrix of $\widetilde{\beta_n}(z)(\alpha_1)$ in the basis $\mathcal{W}$ by checking that \\

$Q(z)(\widetilde{\beta_n}(z)(\alpha_1)-I)_{[\mathcal{V}]}=(\widetilde{\beta_n}(z)(\alpha_1)-I)_{[\mathcal{W}]}Q(z)=$\\

$=z^{n-2} \left( \begin{array}{cc|cc}
	2(z+1)&-(z+1)&&\\
	&&&{\bf 0}\\
	-2&1&&\\
	\hline
	&&&\\
	{\bf 0}&&&{\bf 0}\\
	&&&

\end{array} \right).$\\

	Since ${\widetilde{\beta_n} (z)|}_{B_n}={\widehat{\beta_n} (z)|}_{B_n}=\beta_n(z)$ is  irreducible, and $\widetilde{\beta_n} (z)\neq\widehat{\beta_n} (z)$,     the representations $\widetilde{\beta_n} (z)$ and $\widehat{\beta_n} (z)$ are not equivalent by Lemma 2.9.\\
	
\end{proof}

{\bf Remark 3.3.} One can see from the proof of part 4) of the above lemma the reason why we chose to write the representation $\widetilde{\beta_n} (z)$ in a basis that is different from the one in which  $\widehat{\beta_n} (z)$ is written. Similarly, if one chooses to write  $\widehat{\beta_n} (z)$ in the same basis as $\widetilde{\beta_n} (z),$ one can see that the expression for $\widehat{\beta_n} (z)(\alpha_1)$ is given by

$(\widehat{\beta_n} (z)(\alpha_1))_{[\mathcal{V}]}=I_{n-1}+$\\
\vskip 0.3cm
$+\frac{1}{P_{n}(t)}\left( \begin{array}{cccccc}
	-(1+t)(P_{n-1}(t)+t^{n-2})&&P_{n-1}(t)+t^{n-2}&0&\dots&0\\
	&&&&&\\
	-(1+t)(t^{n-2}-1)&&t^{n-2}-1&0&\dots&0\\
	&&&&&\\
	-(1+t)(t^{n-2}-t)&&t^{n-2}-t&0&\dots&0\\
	&&&&&\\
	\vdots&&\vdots&&\vdots&\\
	&&&&&\\
	-(1+t)(t^{n-2}-t^{n-4})&&t^{n-2}-t^{n-4}&0&\dots&0\\
	&&&&&\\
	-(1+t)(t^{n-2}-t^{n-3})&&t^{n-2}-t^{n-3}&0&\dots&0
\end{array} \right),$  \\  

\vskip .5cm
Thus, it is more convenient to write these representations in two different bases.\\

The existence of the above representations shows that the reduced Burau representation can be extended to the welded braid group. In this section we will show (Theorem 3.25) that for $n$ large enough every extension of the {\it irreducible} reduced Burau representation to $WB_n$ is equivalent to a tensor product of the one-dimensional representation $X_n(1,k),\ k\in \{-1,1\},$   with either $\widetilde{\beta_n} (z)$ or $\widehat{\beta_n} (z),$ where $P_n(z)\neq 0.$\\

Throughout  this section we will work under the following assumptions and use the notations below.

Let $\rho:WB_n\to GL_{n-1}({\mathbb C})$ be an $(n-1)-$dimensional representation of the group $WB_n,$ $n\geqslant 3,$ such that its restriction onto $B_n$ is equivalent to the specialization of the reduced Burau representation $\beta_n(z):B_n\to GL_{n-1}(\mathbb{C})$ for some value $z\in {\mathbb C}^*.$ Suppose that $\beta_n(z)$ is {\it irreducible}, that is $P_n(z)=1+z+z^2+\dots +z^{n-1}\neq 0.$ 

We will select a basis $\mathcal{W}=\{w_1, \dots ,w_{n-1}\}$ of $V,$ such that in this basis $\rho|_{B_n}=\beta_n(z).$ Additionally, we will select a basis $\mathcal{V}=\{v_1, \dots ,v_{n-1}\}$ of $V,$ as defined in Chapter 2 and Lemma 2.6. It remains to determine the action of the generators $\alpha_i$ of the subgroup $S_n$ of $WB_n.$ \\

In the following lemma (Lemma 3.4) we will list some properties, which we will use throughout this section. Each of these statements can be easily obtained from the explicit expression of $\beta_n(z).$ \\

{\bf Lemma 3.4.} Let $\rho:WB_n\to GL_{n-1}({\mathbb C})$ be an $(n-1)-$dimensional representation of the group $WB_n,$ $n\geqslant 3,$ such that its restriction onto $B_n$ is equivalent to the specialization of the reduced Burau representation $\beta_n(z):B_n\to GL_{n-1}(\mathbb{C})$ for some value $z\in {\mathbb C}^*,$ where $P_n(z)\neq 0.$ Then:\\
(a) $dim(Im(C_i))=1,$ $Im(C_i)=span\{v_i\},$ $i=1, \dots, n-1;$\\
(b) $Im(C_{i})\cap Im(C_j)=\{0\}$ for $i\neq j, \  i,j =1,2,\dots ,n-1;$\\
(c) $dim(Ker(C_i))=n-2,$ \\ $Ker(C_i)=span\{w_1,\dots,w_{i-1},w_{i+1},\dots, w_{n-1}\},$ $i=1, \dots, n-1;$\\
(d)  $C_iv_j=0$ for $i,j=1, 2, \dots, n-1,$ $|i-j|\geqslant 2;$\\
(e) $T^n=z^n I.$\\

For the reference, we will supply here the explicit expressions for the operator $T=\rho(\sigma_1\sigma_2\dots \sigma_{n-1})$ and its inverse written in both bases. It is easy to get these matrices by looking at  the consecutive  action of the group elements $\sigma_1,\sigma_2,\dots ,\sigma_{n-1}$ (or their inverses) on the basis vectors, taken in the appropriate order, rather than performing actual matrix multiplication. One can also observe the reason for our choice of the original basis $\mathcal{W}$ in the definition of the reduced Burau representation (see Definition 2.3 and Remark 2.4).\\

{\bf Lemma 3.5.}\\

\noindent
$(T)_{[\mathcal{V}]}=z\cdot\left( \begin{array}{ccccc}
	0&0&\dots&0&-1\\
	1&0&\dots&0&-1\\
	0&1&\dots&0&-1\\
	
	\vdots&\vdots&\ddots&\vdots&\vdots\\
	0&0&\dots&0&-1\\
	0&0&\dots&1&-1
	
\end{array} \right),\ (T^{-1})_{[\mathcal{V}]}=\frac{1}{z}\cdot\left( \begin{array}{ccccc}
	-1&1&\dots&0&0\\
	-1&0&\dots&0&0\\

	\vdots&\vdots&\ddots&\vdots&\vdots\\
	-1&0&\dots&1&0\\
	-1&0&\dots&0&1\\
	-1&0&\dots&0&0
	
\end{array} \right),$
\vskip .3cm

$(T)_{[\mathcal{W}]}=\left( \begin{array}{ccccc}
	-z&-z^2&\dots&-z^{n-2}&-z^{n-1}\\
	1&0&\dots&0&0\\
	0&1&\dots&0&0\\
	
	\vdots&\vdots&\ddots&\vdots&\vdots\\
	0&0&\dots&0&0\\
	0&0&\dots&1&0
	
\end{array} \right),$

\vskip .3cm

$(T^{-1})_{[\mathcal{W}]}=\left( \begin{array}{cccrr}
	0&1&\dots&0&0\\
	0&0&\dots&0&0\\

	\vdots&\vdots&\ddots&\vdots&\vdots\\
	0&0&\dots&1&0\\
	0&0&\dots&0&1\\
	-\frac{1}{z^{n-1}}&	-\frac{1}{z^{n-2}}&\dots&	-\frac{1}{z^{2}}&	-\frac{1}{z}
	
\end{array} \right).$
\vskip .3cm

{\bf Lemma 3.6.} For $n\geqslant 4,$

1) Each of the vectors $v_3, v_4, \dots, v_{n-1}$ is an eigenvector of $A_1.$

2) All vectors $ v_3, v_4, \dots, v_{n-1}$  have the same  eigenvalue $k,$ where $k=1$ or $k=-1.$\\

\begin{proof} 1)  For each $j=3, \dots ,n-1,$ by Lemma 3.4(a), \\ $v_j\in Im(C_{j}).$ Thus, there exists a vector $x_j,$ such that $v_j=C_jx_j.$
	
	By Lemma 2.11(b), $A_1v_j=A_1C_jx_j=C_jA_1x_j\in Im(C_{j})=span\{v_{j}\}.$\\
	
	2) We will show that if the set $\{ v_3, v_4, \dots, v_{n-1}\}$ consists of more than one vector ($n\geqslant 5$), then any two consecutive vectors have the same eigenvalue.\\
	
	Let $A_1v_j=k_jv_j$ for $j=3,4,\dots ,n-1,$ where $k_j\in  \{-1,1\}$ (Lemma 2.11(f)).
	
	For every $j=3, 4, \dots , n-2$ consider $A_2v_{j+1}.$ By using Lemmas 2.11(e) and 3.5, we have\\
	$A_2v_{j+1}=TA_1T^{-1}v_{j+1}=\frac{1}{z}T(A_1v_j)=\frac{1}{z}\cdot k_jTv_j=\frac{1}{z}\cdot k_j\cdot zv_{j+1}=k_jv_{j+1}.$\\
	
	Rewriting  the group identity (3) as $\alpha_i =\alpha^{-1}_{i+1}\alpha^{-1}_i\alpha_{i+1}\alpha_i \alpha_{i+1}$  we get:\\
	
	$k_{j+1}v_{j+1}=A_1v_{j+1}=A_2^{-1}A_1^{-1}A_2A_1A_2v_{j+1}=k_j^{-1}k_{j+1}^{-1}k_jk_{j+1}k_jv_{j+1}=$
	
	$=k_jv_{j+1},$ and since $v_{j+1}\neq 0$ it follows that $k_j=k_{j+1}.$\\

\end{proof}

{\bf Corollary 3.8.} For $n\geqslant 4,$ $A_iv_j=kv_j$ for all $i,j=1, 2, \dots, n-1,$ $|i-j|\geqslant 2,$ where $k\in \{-1,1\}.$\\

\begin{proof}  By Lemmas 3.5 and 2.11(g), for every $i\neq j, \ i,j=1, \dots, n-1,$\\
	$(T^{-1})^{i}v_j\in span\{v_{[j-i]}\}$, where $[m]=m\pmod n.$ 
	
	Let $k$ be the eigenvalue of $A_1$ for the eigenvectors $v_3, \dots , v_{n-1}.$ Then for all $i,j=1, \dots , n-1,$ $|i-j|\geqslant 2,$  \\
	
	$(A_i-kI)v_j=T^{i-1}(A_1-kI)(T^{-1})^{i-1}v_j=0,$ since \\$(T^{-1})^{i-1}v_j\in span \{v_{[j-i+1]}\}\,$ and  $[j-i+1]\neq 0,1,2.$\\

\end{proof}

{\bf Corollary 3.8.} For $n\geqslant 4,$ $(A_i-kI)C_j=C_j(A_i-kI)=0$ for all $i,j=1, 2, \dots, n-1,$ $|i-j|\geqslant 2,$ where $k\in \{-1,1\}.$\\

\begin{proof}
	By Lemma 3.4(a), $Im(C_j)=span\{v_j\},$ so $(A_i-kI)C_j=0$ for all $i,j=1, 2, \dots, n-1,$ $|i-j|\geqslant 2,$ and by Lemma 2.11(b), \\$(A_i-kI)C_j=C_j(A_i-kI).$\\
\end{proof}

{\bf Lemma 3.9.} Let $\rho:WB_n\to GL_{n-1}({\mathbb C})$ be an $(n-1)-$dimensional representation of the group $WB_n$ for $n\geqslant 5,$ such that its restriction onto $B_n$ is equivalent to the specialization of the reduced Burau representation $\beta_n(z):B_n\to GL_{n-1}(\mathbb{C})$ for some value $z\in {\mathbb C}^*$ (not necessarily irreducible). Then the restriction of $\rho$ onto $S_n$ is {\it irreducible}.

\begin{proof} Suppose ${\rho|}_{S_n}$ is reducible. By the theorem of Burnside (\cite{Burnside}, note C, p. 468), for $n\geqslant 5,$ the symmetric group $S_n$ has no irreducible representations of dimension $r$ for $2\leqslant r\leqslant n-2.$ Thus, by Maschke's theorem,  ${\rho|}_{S_n}$ is a direct sum of one-dimensional representations. Every one-dimensional representation $\eta$ of $S_n$  is equivalent to either identity representation $id_n:S_n\to 1,$ or a sign representation $sgn_n:S_n\to \{-1,1\}.$ In either case  $\eta(\alpha_i)=\eta(\alpha_j)$ for all $1\leqslant i,j\leqslant n-1$. Hence, $\rho(\alpha_i)=\rho(\alpha_j)$ for all $1\leqslant i,j\leqslant n-1.$
	Then, from the group identities (7) and (1), it follows that $\rho(\sigma_i)=\rho(\sigma_{i+1})$,  for $i=1, 2,\dots,n-2,$ a contradiction with the fact that ${\rho|}_{B_n}$ is equivalent to the specialization of the reduced Burau representation.\\
	
\end{proof}

Denote by $\xi_n:S_n\to GL_{n-1}(\mathbb{C})$  the $(n-1)-$ dimensional irreducible representation of $S_n,$ which is an $(n-1)-$dimensional irreducible subrepresentation of the $n-$dimensional natural permutation representation on the subspace of those vectors whose coordinates add up to zero. (This representation is sometimes called the standard representation of the symmetric group; it is not to be mixed up with the $n-$dimensional standard representation of the braid group $\tau_n(t).$)\\

{\bf Corollary 3.10.} Suppose $n=5$ or  $n\geqslant 7.$ \\Let $\rho:WB_n\to GL_{n-1}({\mathbb C})$ be an $(n-1)-$dimensional representation of the group $WB_n,$ for  such that ${\rho|}_{B_n}$ is equivalent to the specialization of the reduced Burau representation $\beta_n(z):B_n\to GL_{n-1}(\mathbb{C})$ for some value $z\in {\mathbb C}^*.$ Then ${\rho|}_{S_n}$ is equivalent to either $(n-1)-$ dimensional irreducible representation  $\xi_n:S_n\to GL_{n-1}(\mathbb{C}),$  or its tensor product with the  one-dimensional sign representation, $sgn_n\otimes \xi_n:S_n\to GL_{n-1}(\mathbb{C}).$\\

\begin{proof} This statement immediately follows from Lemma 3.9 and a well-known fact from the theory of representations of symmetric groups (see, for example, \cite{JK}, Theorem 2.4.10(ii)), that for $n\geqslant 5,$   $n\neq 6, $ the listed  representations are the only two  non-equivalent irreducible $(n-1)-$dimensional representations of $S_n.$ 
	
\end{proof}

First, we are going to consider only those representations $\rho,$ whose restriction to $S_n$ is equivalent to $\xi_n$ (and not to its tensor product with  the  sign representation), which are exactly the representations (as Lemma 3.11 below shows) with $A_1v_j=v_j, \ j=3, 4, \dots, n-1.$  The complete classification will follow then from the classification theorem for this specific case.\\

{\bf Lemma 3.11.} Let $\rho:WB_n\to GL_{n-1}({\mathbb C})$ be an $(n-1)-$dimensional representation of the group $WB_n$ for  $n=5$ or  $n\geqslant 7,$ such that ${\rho|}_{B_n}$ is equivalent to the specialization of the reduced Burau representation $\beta_n(z):B_n\to GL_{n-1}(\mathbb{C})$ for some value $z\in {\mathbb C}^*,$ such that  $P_n(z)\neq 0.$ Then ${\rho|}_{S_n}$ is equivalent to  $\xi_n:S_n\to GL_{n-1}(\mathbb{C})$ if and only if the vectors $v_3, v_4, \dots, v_{n-1}$ are the eigenvectors of $\rho(\alpha_1)$ with the same eigenvalue 1.\\

\begin{proof} The eigenvalues of $\xi_n(\alpha_1)$ are 1 of multiplicity (n-2), and -1 of multiplicity 1. By Lemma 3.6, $\rho(\alpha_1)$ has eigenvalue $k$ of multiplicity at least 2, hence $k=1.$
	
\end{proof}

Consider now the operators \\$D_i=\rho(\alpha_i)-I=A_i-I,$ for $ i=1, 2, \dots, n-1,$\\ introduced in Chapter 2.\\

{\bf Lemma 3.12.} Suppose $\rho:WB_n\to GL_{n-1}({\mathbb C})$ for $n=5$ or $n\geqslant 7,$ such that  $\rho|_{B_n}$ is equivalent to $\beta_n(z),$ where $z\in \mathbb{C^*},$  $P_n(z)\neq 0.$ Suppose that  $\rho(\alpha_1)v_j=v_j$ for $j=3, 4, \dots , n-1.$ \\
Then for all $i,j=1, 2, \dots, n-1$ :\\
(a) $dim(Im(D_i))=1,$  $dim(Ker(D_i))=n-2;$\\
(b) $Im(D_i)\cap Im(D_j)=\{0\},$ $i\neq j;$ \\
(c) $D_iD_j=0,$ $|i-j|\geqslant 2;$\\
(d) $D_i^2=-2D_i;$\\
(e) $D_iD_jD_i-D_i=0,$ $|i-j|=1;$\\
(f)  $D_iC_j=C_jD_i=0,$ $|i-j|\geqslant 2.$\\

\begin{proof}  By Lemma 3.11, ${\rho|}_{S_n}$ is equivalent to  $\xi_n.$ The representation $\xi_n,$ when precomposed with the natural homomorphism $B_n\to S_n,$    is equivalent  to the specialization of Burau representation $\beta_n(1)$.   Now each of the statements (a)--(e) follow easily from the explicit expressions for the matrices of $\beta_n(1)$. The statement (f) follows from Corollary 3.8.  \\
	
\end{proof}

{\bf Lemma 3.13.} Under conditions of Lemma 3.12, \\$D_1v_2\neq 0.$

\begin{proof} Suppose that $D_1v_2=0.$ Then $A_1v_2=v_2,$ and, by Lemmas 2.11(e) and 3.5, 
	$A_2v_3=\frac{1}{z}T(A_1v_2)=\frac{1}{z}T(v_2)=v_3.$\\
	
	By Lemma 2.11(d), $C_1A_{2}A_{1}= A_{2}  A_{1}C_{2},$ and by Lemma 3.4(d), \\$C_1v_3=0.$  Since $C_2v_3=v_2$ and $A_1v_3=v_3,$ we obtain:\\
	
	$0=(A_{2}  A_{1}C_{2}-C_1A_{2}A_{1})v_3=A_{2}  A_{1}v_2-C_1A_{2}v_3=A_{2} v_2-C_1v_3=$\\
	$=A_{2} v_2\neq 0,$ a contradiction.
	
\end{proof}

{\bf Lemma 3.14.} Under conditions of Lemma 3.12, \\$Im(D_1)\subseteq \bigcap\limits_{j=3}^{n-1} Ker(C_j).$

\begin{proof} By Lemma 3.12(f), for every  $j=3, 4, \dots, n-1,$ \\we have $C_jD_1=0,$ hence $Im(D_1)\subseteq Ker(C_j).$\\
	
\end{proof}

{\bf Lemma 3.15.} $dim(\bigcap\limits_{j=3}^{n-1} Ker(C_j))=2, $ and \\$\bigcap\limits_{j=3}^{n-1} Ker(C_j)=span \{w_1, w_2\}.$\\

\begin{proof}
	Immediately follows from the facts that $dim(Ker(C_i))=n-2,$ and $Ker(C_i)=span\{w_1,\dots,w_{i-1},w_{i+1},\dots, w_{n-1}\},$ $i=1, \dots, n-1$ (Lemma 3.4(c)).\\
\end{proof}

Now we will write the basis of $\bigcap\limits_{j=3}^{n-1} Ker(C_j)$ in terms of vectors \\$v_1, \dots , v_{n-1}.$\\

{\bf Lemma 3.16.}   $\bigcap\limits_{j=3}^{n-1} Ker(C_j)=span\{v_1,{\bf v}\},$ where \\${\bf v}=(1+z+\dots +z^{n-3})v_2+z (1+z+\dots +z^{n-4})v_3+\dots +$\\$+z^{n-4}(1+z)v_{n-2}+z^{n-3}v_{n-1}=\sum\limits_{i=2}^{n-1}z^{i-2}P_{n-i}(z)v_{i}.$

\begin{proof} Rather than using the change-of-basis matrix $Q(z),$ instead we can solve the system of $(n-3)$ linear equations (of rank $(n-3)$) in $(n-1)$ variables, using the explicit expressions for the matrices $C_j:$ the vector $v=\sum\limits_{i=1}^{n-1}x_iv_i \in Ker(C_j)\Longleftrightarrow zx_{j-1}+(-z-1)x_j+x_{j+1}=0$ for every $j=3, 4,\dots, n-2,$ and \\ $v=\sum\limits_{i=1}^{n-1}x_iv_i \in Ker(C_{n-1})\Longleftrightarrow zx_{n-2}+(-z-1)x_{n-1}=0.$ \\
	It is easy to see that $C_jv_1=0$ and $C_j{\bf v}=0$ for every $j=3, 4,\dots, n-1,$ thus two linearly independent vectors $v_1$ and ${\bf v}$ form a basis of two-dimensional space $\bigcap\limits_{j=3}^{n-1} Ker(C_j).$ \\ 
\end{proof}

For the one-dimensional subspaces $Im(C_1)$ and $Im(D_1),$ the following two cases are possible:\\

\begin{flushleft}

	{\bf Case I:} $Im(C_1)=Im(D_1).$ \\
	{\bf Case II:} $Im(C_1)\neq Im(D_1).$ \\
\end{flushleft}

We will consider these two cases separately. We will show that, under conditions of Lemma 3.12, in Case II  $Ker(C_1)=Ker(D_1)$ (see Lemma 3.21). Moreover, we will prove that in each of these cases for every $z\neq 1, \ P_n(z)\neq 0,$ $ z \in \mathbb{C}^*,$ there exists exactly one extension of the irreducible reduced Burau representation $\beta_n(z):$  $\widetilde{\beta_n} (z)$ and $\widehat{\beta_n} (z)$ respectively; and for $z=1,$ there exists exactly one extension $\widetilde{\beta_n} (1)=\widehat{\beta_n} (1)$ for which $Im(C_1)=Im(D_1)$ and $Ker(C_1)=Ker(D_1).$\\

\vskip 0.5cm

{\bf Case I.} \\ 

{\bf Lemma 3.17.} Suppose $\rho:WB_n\to GL_{n-1}({\mathbb C})$ for $n=5$ or $n\geqslant 7,$ such that  $\rho|_{B_n}$ is equivalent to $\beta_n(z),$ where $z\in \mathbb{C^*},$  $P_n(z)\neq 0.$ Suppose that  $\rho(\alpha_1)v_j=v_j$ for $j=3, 4, \dots , n-1.$ Suppose $Im(D_1)=Im(C_1).$ \\

Then $D_1v_1=-2v_1,$ and $D_1v_2=v_1.$\\

\begin{proof}  
	
	Since $v_1\in Im(C_1)=Im(D_1),$ there exists a vector $v',$ such that $v_1=D_1v'.$ By Lemma 3.12(d),\\ $D_1v_1=D_1^2v'=-2D_1v'=-2v_1.$ \\

	The non-zero vector $D_1v_2\in Im(D_1)=Im(C_1)=span\{v_1\},$ hence
	$D_1v_2=\lambda v_1,$ where $\lambda=\lambda(z)\in\mathbb{C^*}.$\\

	Consider vector $D_2v_1.$ By Lemmas 3.12(d) and 3.5,\\
	
	$D_2v_1=TD_1T^{-1}v_1=-\frac{1}{z}\cdot TD_1(v_1+v_2+\dots +v_{n-1})=$\\
	$=-\frac{1}{z}\cdot T(D_1v_1+D_1v_2)=-\frac{1}{z}\cdot T(-2v_1+\lambda v_1)=$\\
	$=\frac{2-\lambda}{z}Tv_1=\frac{2-\lambda}{z}(zv_2)=(2-\lambda)v_2.$ \\
	
	By Lemma 3.12(e), $D_1D_2v_1=D_1D_2D_1v'=D_1v'=v_1,$ so\\

	$0=D_1D_2v_1-v_1=D_1(2-\lambda)v_2-v_1=(\lambda(2-\lambda)-1)v_1,$ hence \\
	$\lambda(2-\lambda)-1=0$ (since $v_1\neq 0$).\\  Thus $\lambda=1,$ and $D_1v_2=v_1.$ \\
	
\end{proof}

{\bf Lemma 3.18.} Suppose $\rho:WB_n\to GL_{n-1}({\mathbb C})$ for $n=5$ or $n\geqslant 7,$ such that  $\rho|_{B_n}$ is equivalent to $\beta_n(z),$ where $z\in \mathbb{C^*},$  $P_n(z)\neq 0.$ Suppose that  $\rho(\alpha_1)v_j=v_j$ for $j=3, 4, \dots , n-1.$ Suppose $Im(D_1)=Im(C_1).$ \\

Then $\rho$ is equivalent to $\widetilde{\beta_n} (z).$\\

\begin{proof} 	By Lemma 3.17, \\

	$(\rho(\alpha_1))_{[\mathcal{V}]}=(I+D_1)_{[\mathcal{V}]}=\left( \begin{array}{rcc|c}
		-1&&1&\\
		&&&0\\
		0&&1&\\
		\hline
		&&&\\
		&0&&I_{n-3}\\
		&&&
	\end{array} \right)=\widetilde{\beta_n} (z)(\alpha_1).$ \\ 
	\vskip 0.5cm 
	
By Lemma 2.11(e), definition of $\widetilde{\beta_n} (z)$ (Definition 3.1) and the explicit expressions for the matrices $T$ and $T^{-1}$ in basis $\mathcal{V}$ (Lemma 3.5) we obtain $(\rho(\alpha_i))_{[\mathcal{V}]}=	\widetilde{\beta_n} (z)(\alpha_i)$ for all $i=1, 2,\dots, n-1,$ and hence, ${(\rho)_{[\mathcal{V}]}}=\widetilde{\beta_n} (z).$\\

\end{proof}

{\bf Case II.} \\  

Consider a non-zero vector $D_1v_2$   spanning one-dimensional subspace $Im(D_1).$ Since in this case $Im(D_1)\neq Im(C_1)=span\{v_1\},$ it follows that $D_1v_2\notin span\{v_1\}.$ 	By Lemmas 3.14 and 3.16, \\$D_1v_2\in Im(D_1)\subseteq \bigcap\limits_{j=3}^{n-1} Ker(C_j)=span\{v_1,{\bf v}\},$ where \\${\bf v}=\sum\limits_{i=2}^{n-1}z^{i-2}P_{n-i}(z)v_{i}.$ We can write $D_1v_2=av_1+b\sum\limits_{i=2}^{n-1}z^{i-2}P_{n-i}(z)v_{i},$ where $a=a(z), b=b(z)\in \mathbb{C}.$ Moreover, since $D_1v_2\notin span\{v_1\},$ it follows that $b\neq 0.$\\

{\bf Lemma 3.19.} Suppose $\rho:WB_n\to GL_{n-1}({\mathbb C})$ for $n=5$ or $n\geqslant 7,$ such that  $\rho|_{B_n}$ is equivalent to $\beta_n(z),$ where $z\in \mathbb{C^*},$  $P_n(z)\neq 0.$ Suppose that  $\rho(\alpha_1)v_j=v_j$ for $j=3, 4, \dots , n-1.$ Suppose $Im(D_1)\neq Im(C_1).$  \\ Then $D_1v_1=-(1+z)D_1v_2.$\\

\begin{proof}

	By Corollary 3.7, $D_{n-1}v_j=0$ for all $j=1, 2,\dots, n-3,$
	and by Lemma 3.12(c), $D_{n-1}D_1=0.$ Thus, using Lemma 2.11(e),  \\
	
	\noindent
	$0=D_{n-1}D_1v_2=D_{n-1}\left[av_1+b \sum\limits_{i=2}^{n-1}z^{i-2}P_{n-i}(z)v_{i}\right]=$\\ $=bD_{n-1}(z^{n-4}P_2(z)v_{n-2}+z^{n-3}P_1(z)v_{n-1})=$\\
	$=bz^{n-4}D_{n-1}((1+z)v_{n-2}+zv_{n-1})=$\\$=bz^{n-4}T^{n-2}D_{1}\left( T^{-1}\right) ^{n-2}\left((1+z)v_{n-2}+zv_{n-1}\right).$ \\Since $b\neq 0$ and $z\neq 0,$ and $T$ is invertible, we obtain\\
	
	$D_{1}\left((1+z)\left( T^{-1}\right) ^{n-2}v_{n-2}+z\left( T^{-1}\right) ^{n-2}v_{n-1}\right)=0.$\\
	
	By Lemma 3.5, \\$\left( T^{-1}\right) ^{n-2}v_{n-2}=\frac{1}{z^{n-3}}\left( T^{-1}\right) v_1=-\frac{1}{z^{n-2}}(v_1+v_2+\dots v_{n-1}),$\\
	so by Corollary 3.7, $D_1\left( T^{-1}\right) ^{n-2}v_{n-2}=-\frac{1}{z^{n-2}}(D_1v_1+D_1v_2).$\\
	Again, by Lemma 3.5, $\left( T^{-1}\right) ^{n-2}v_{n-1}=\frac{1}{z^{n-2}}v_1.$ Thus, \\
	$-(1+z)(D_1v_1+D_1v_2)+zD_1v_1=0.$\\
	Solving for $D_1v_1,$ we obtain $D_1v_1=-(1+z)D_1v_2.$\\
	
\end{proof}	

{\bf Corollary 3.20.} Under the conditions of Lemma 3.19, \\$v_1+(z+1)v_2\in Ker(D_1).$\\

{\bf Lemma 3.21.} Under the conditions of Lemma 3.19, \\$Ker(D_1)=Ker(C_1).$\\

\begin{proof} By Lemma 3.12(a), $dim(Ker(D_1))=n-2.$ We claim that $n-2$ linearly independent vectors $v_1+(z+1)v_2, v_3, v_4, \dots , v_{n-1}$ form a basis of $Ker(D_1).$ Indeed, by the conditions of the lemma, $v_j\in Ker(D_1)$ for \\$j=3, 4, \dots , n-1,$ and by Corollary 3.20, $v_1+(z+1)v_2\in Ker(D_1).$ By Lemma 3.4(d), 
	$v_j\in Ker(C_1)$ for $j=3, \dots , n-1,$ and by the explicit expression for $(C_1)_{[\mathcal{V}]},$\\ $C_1(v_1+(z+1)v_2)=-(z+1)v_1+(z+1)v_1=0.$ Thus, \\$Ker(D_1)\subseteq Ker(C_1), $ and since\\ $dim(Ker(D_1))=dim(Ker(C_1))=n-2, $ we have $Ker(D_1)=Ker(C_1).$ \\
	
\end{proof}

{\bf Corollary 3.22.} Under the conditions of Lemma 3.19, \\$Ker(D_1)=span\{w_2,w_{3},\dots, w_{n-1}\}. $\\

\begin{proof}
	Immediately follows from Lemma 3.21 and Lemma 3.4(c).
\end{proof}

{\bf Lemma 3.23.}  Suppose $\rho:WB_n\to GL_{n-1}({\mathbb C})$ for $n=5$ or $n\geqslant 7,$ such that  $\rho|_{B_n}$ is equivalent to $\beta_n(z),$ where $z\in \mathbb{C^*},$  $P_n(z)\neq 0.$ Suppose that  $\rho(\alpha_1)v_j=v_j$ for $j=3, 4, \dots , n-1.$ Suppose $Im(D_1)\neq Im(C_1).$ \\
1) If $z\neq 1,$ then $D_1w_1=-2w_1+\frac{1}{z} w_2;$\\
2) If $z=1,$ then no such representation $\rho$ exists.\\

\begin{proof} By Corollary 3.22, $D_1w_1\neq 0.$ By Lemmas 3.14 and 3.15, \\$Im(D_1)\subseteq \bigcap\limits_{j=3}^{n-1} Ker(C_j)=span\{w_1,w_2\}.$ Thus, we can write \\$D_1w_1=a_1w_1+a_2w_2,$ where $a_1=a_1(z),a_2=a_2(z)\in \mathbb{C},$ with at least one of  $a_1,a_2$ not equal to zero. Additionally, since\\ $Im(D_1)\neq Im(C_1)= span \{v_1\},$ \\$D_1w_1=a_1w_1+a_2w_2\notin span \{v_1\}$ and $v_1= -z^{n-2}(z+1)w_1+z^{n-2}w_2,$ we have $a_1+(z+1)a_2\neq 0.$\\

	By Lemma 3.12(d) and Corollary 3.22,\\ $-2D_1w_1=D_1^2w_1=D_1(a_1w_1+a_2w_2)=a_1D_1w_1,$ and since $D_1w_1\neq 0,$ we have $a_1=-2.$ So, \\ $D_1w_1=-2w_1+a_2w_2,$ where $-2+(z+1)a_2\neq 0.$\\

	By Lemma 3.12(e), \\
	$D_1D_2D_1-D_1=0,$ so, using Lemmas 3.5, 2.11(e), and Corollary 3.22, \\
	$D_1w_1=D_1D_2D_1w_1=D_1D_2(-2w_1+a_2w_2)=\\ =D_1TD_1(-2T^{-1}w_1+a_2T^{-1}w_2)=\\ =D_1TD_1\left(\frac{2}{z^{n-1}}w_{n-1}+a_2(w_1-\frac{1}{z^{n-2}}w_{n-1})\right) =\\ =a_2D_1TD_1w_1=a_2D_1T(-2w_1+a_2w_2)=\\ =-2a_2D_1(-zw_1+w_2)+a_2^2D_1(-z^2w_1+w_3)=\\=(2a_2z-a_2^2z^2)D_1w_1.$\\
	
	Again, since $D_1w_1\neq 0,$ we obtain \\$1-2a_2z+a_2^2z^2=0,$ so \\$(a_2z-1)^2=0,$ hence \\ $a_2=\frac{1}{z}.$\\
	
	Since for $Im(D_1)\neq Im(C_1)$ we must have  that $-2+(z+1)\cdot \frac{1}{z}\neq 0,$ it follows that $z\neq 1.$ If $z=1,$ then $D_1w_1=-2w_1+w_2=v_1,$ so no representation for which $Im(D_1)\neq Im(C_1)$ exists.\\

\end{proof}

{\bf Lemma 3.24.} Suppose $\rho:WB_n\to GL_{n-1}({\mathbb C})$ for $n=5$ or $n\geqslant 7,$ such that  $\rho|_{B_n}$ is equivalent to $\beta_n(z),$ where $z\in \mathbb{C^*},$  $P_n(z)\neq 0,$ and $z\neq 1.$ Suppose that  $\rho(\alpha_1)v_j=v_j$ for $j=3, 4, \dots , n-1.$ Suppose $Im(D_1)\neq Im(C_1).$ \\

Then $\rho$ is equivalent to $\widehat{\beta_n} (z).$\\

\begin{proof} By Lemma 3.23, \\

	$(\rho(\alpha_1))_{[\mathcal{W}]}=(I+D_1)_{[\mathcal{W}]}=\left( \begin{array}{rcc|c}
		-1&&0&\\
		&&&0\\
		\frac{1}{z}&&1&\\
		\hline
		&&&\\
		&0&&I_{n-3}\\
		&&&
	\end{array} \right)=\widehat{\beta_n} (z)(\alpha_1).$ \\  
	\vskip 0.5cm 
	
	Like in Lemma 3.18, by using Lemma 2.11(e), definition of $\widehat{\beta_n} (z)$ (Definition 3.1) and the explicit expressions for the matrices $T$ and $T^{-1}$ in basis $\mathcal{W}$ (Lemma 3.5) we obtain $(\rho(\alpha_i))_{[\mathcal{W}]}=	\widehat{\beta_n} (z)(\alpha_i)$ for all $i=1, 2,\dots, n-1,$ and hence, ${(\rho)_{[\mathcal{W}]}}=\widehat{\beta_n} (z).$\\

\end{proof}

{\bf Theorem 3.25.} Let $n=5$ or $n\geqslant 7.$ Suppose $\rho:WB_n\to GL_{n-1}({\mathbb C})$ is an $(n-1)-$dimensional complex representation of welded braid group $WB_n,$  such that its restriction onto braid group $\rho|_{B_n}$ is equivalent to a specialization of a reduced Burau representation $\beta_n(z)$ for some $z \in\mathbb{C}^*,$ where $P_n(z)\neq 0.$\\

Then there exists a unique pair $(k,\beta'),$ such that  $\rho$ is equivalent to $X_n(1,k)\otimes\beta',$ where $k\in \{ -1,1\},$ and $\beta'$ is a representation of $WB_n,$\\  $\beta'\in \{\widetilde{\beta_n} (z),\widehat{\beta_n} (z)\},$ if $z\neq 1,$ and \\$\beta'=\widetilde{\beta_n} (1),$ if $z=1.$\\

\begin{proof} Select a basis  $\mathcal{W}=\{w_1, \dots ,w_{n-1}\}$ of $(n-1)-$dimensional complex space  $V,$ such that in this basis $\rho|_{B_n}=\beta_n(z).$
	Additionally, select another basis  $\mathcal{V}=\{v_1, \dots ,v_{n-1}\}$ using the change-of-basis matrix from Lemma 2.6.\\
	
	By Lemma 3.6, there exists $k\in \{ -1,1\},$ such that $\rho(\alpha_1)v_j=kv_j$ for $j=3, 4, \dots , n-1.$\\
	
	Consider $\hat{\rho}=(X_n(1,k))^{-1}\otimes {\rho},$ where $X_n(1,k)$ is a one-dimensional representation of $WB_n$ (see Definition 2.7). Then $\hat{\rho}$ is an \\$(n-1)-$dimensional representation of $WB_n,$ such that $\hat{\rho}_{|B_n}$ is equivalent to $\beta_n(z)$ and \\$\hat{\rho}(\alpha_1)v_j=v_j$ for all $j=3,4,\dots, n-1.$\\
	
	Consider $D_1=\hat{\rho}(\alpha_1)-I,$ and $C_1=\hat{\rho}(\sigma_1)-I.$ By Lemma 3.12(a), $dim(Im(D_1))=1,$ and by Lemma 3.4(a), $dim(Im(C_1))=1.$\\
	
	If $Im(D_1)=Im(C_1),$ then by Lemma 3.18, $ \hat{\rho}$ is equivalent to $\widetilde{\beta_n} (z).$ \\
	
	If $Im(D_1)\neq Im(C_1),$ then by Lemma 3.23, $z\neq 1,$ and by Lemma 3.24, $ \hat{\rho}$ is equivalent to $\widehat{\beta_n} (z).$ \\
	
	Thus, $\rho=X_n(1,k)\otimes \hat{\rho}$ is equivalent to $X_n(1,k)\otimes\beta',$ where $\beta'$ is either $\widetilde{\beta_n} (z)$ or $\widehat{\beta_n} (z).$\\
	
	Now we will show that for distinct pairs $(k,\beta')$ the representations $X_n(1,k)\otimes\beta'$ are not equivalent. Indeed, the operator \\$(X_n(1,k)\otimes\beta')(\alpha_1)$ has eigenvalue $k$ of multiplicity $(n-2),$ and eigenvalue $-k$ of multiplicity $1,$ hence for distinct values of $k,$ the pairs $(k,\beta')$ give non-equivalent representations. 
	
	Next, for $z\neq 1,$ by Lemma 3.2, part 4), the representations $\widetilde{\beta_n} (z)$ and $\widehat{\beta_n} (z)$ are not equivalent, so, for every $k\in \{ -1,1\},$  the representations $X_n(1,k)\otimes\widetilde{\beta_n} (z)$ and $X_n(1,k)\otimes\widehat{\beta_n} (z)$ are not equivalent.

\end{proof}

{\bf Remark 3.26.} Note that for $z=1,$ since $\widetilde{\beta_n} (1)$  is equivalent to $ \widehat{\beta_n} (1),$ the corresponding part  of the theorem can also be formulated as $\beta'=\widehat{\beta_n} (1).$

 \section{Extensions of the standard irreducible representation of dimension $n$ for $n\geqslant 4.$}
 
 The goal of the next two sections is to classify all non-equivalent extensions of the $n-$dimensional standard {\it irreducible} representation of the braid group $B_n$ to the welded braid group $WB_n$ for $n\geqslant 3.$ Current section is devoted to the case $n\geqslant 4,$ and Section 5 is devoted to the case $n=3.$

 We will start this section with the following definition.\\
 
 {\bf Definition 4.1.} Let $WB_n$ be the welded braid group, $n\geqslant 3.$ \\ Define 
 $\widetilde{\tau}_n (t,q):WB_n\to GL_n(\mathbb{Z}[t^{\pm 1},q^{\pm 1}]) $ by\\
 
 \noindent
 $\widetilde{\tau}_n (t,q)(\sigma_i)=\left( \begin{array}{c|ccc|c}
 	&&&&\\
 	I_{i-1}&&0&&0\\
 	&&&&\\
 	\hline 
 	
 	&0&&t&\\
 	0&&&&0\\
 	&1&&0&\\
 	\hline
 	&&&&\\
 	0&&0&&I_{n-1-i}\\
 	&&&&
 \end{array} \right),$ 
 \vskip .5cm
 
 \noindent
 $\widetilde{\tau}_n(t,q)(\alpha_i)=\left( \begin{array}{c|ccc|c}
 	&&&&\\
 	I_{i-1}&&0&&0\\
 	&&&&\\
 	\hline 
 	&&&&\\
 	&0&&\frac{1}{q}&\\
 	0&&&&0\\
 	&q&&0&\\
 	&&&&\\
 	\hline
 	&&&&\\
 	0&&0&&I_{n-1-i}\\
 	&&&&
 \end{array} \right),$
 \vskip .5cm
 
 \noindent
 for $i=1,2,\dots, n-1,$ where $t$ and $q$ are indeterminates, and $I_{m}$ is the $m\times m$ identity matrix.\\
 
 For any $z,\lambda \in \mathbb{C}^*,$ denote by 
 $\widetilde{\tau}_n(z,\lambda):WB_n\to GL_n(\mathbb{C})$ 
 the map obtained from $\widetilde{\tau}_n(t,q)$ by specializing $t\mapsto z$ and $q\mapsto\lambda.$ \\
 
 {\bf Lemma 4.2.} 1)  $\widetilde{\tau}_n(t,q):WB_n\to GL_n(\mathbb{Z}[t^{\pm 1},q^{\pm 1}]) $  is a representation of $WB_n,$ which is an extension of the standard representation $\tau_n(t)$  of the braid group  $B_n.$\\
 
 2) For all $z,\lambda \in \mathbb{C}^*,$ the map  $\widetilde{\tau}_n(z,\lambda)$ is an $n-$dimensional complex representation of $WB_n,$ which is an extension of  $\tau_n(z).$  If  $z\neq 1,$ then $\widetilde{\tau}_n(z,\lambda)$ is irreducible.\\
 
 3) For all $z_1,z_2,\lambda_1,\lambda_2 \in \mathbb{C}^*:$\\
 a) if $z_1\neq z_2,$ then $\widetilde{\tau}_n(z_1,\lambda_1)$ and $\widetilde{\tau}_n(z_2,\lambda_2)$ are  not equivalent;\\
 b) if $z_1=z_2=z\neq 1$ and $\lambda_1\neq \lambda_2,$ then $\widetilde{\tau}_n(z,\lambda_1)$ and $\widetilde{\tau}_n(z,\lambda_2)$ are  not equivalent.\\
 
 \begin{proof} 1) It is easy to check that all defining relations (1)--(8) are satisfied, and $\widetilde{\tau}_n(t,q)(\sigma_i)=\tau_n(t)(\sigma_i)\ \ \forall i=1,\dots,n-1.$\\
 	
 2) The first statement immediately follows from part 1).  By Theorem 2.5, part 2), if $z\neq 1,$ then $\tau_n(z)$  is irreducible, hence $\widetilde{\tau}_n(z,\lambda)$ is irreducible.\\
 
 3)  Since $det(\tau_n(z)(\sigma_1))=-z,$ the representations $\tau_n(z_1)$ and $\tau_n(z_2)$ are not equivalent for $z_1\neq z_2,$ hence $\widetilde{\tau}_n(z_1,\lambda_1)$ and $\widetilde{\tau}_n(z_2,\lambda_2)$ are  not equivalent.\\
 For $z_1=z_2=z\neq 1,$ the restriction to the braid group $\tau_n(z)$ is irreducible. Then by Lemma 2.9,  since $\widetilde{\tau}_n(z,\lambda_1)\neq\widetilde{\tau}_n(z,\lambda_2)$ for  $\lambda_1\neq \lambda_2,$ the representations $\widetilde{\tau}_n(z,\lambda_1)$ and $\widetilde{\tau}_n(z,\lambda_2)$ are  not equivalent.

 \end{proof}

 The existence of the above representation shows that the standard representation can be extended to the welded braid group. In this section we will show (Theorem 4.11) that for $n$ large enough every extension of the {\it irreducible} standard representation to $WB_n$ is equivalent to a tensor product of the one-dimensional representation $X_n(1,k),$  $k\in \{-1,1\},$   and $\widetilde{\tau}_n(z,\lambda),$  where $z\neq 1, \ z\in {\mathbb C}^*.$\\
 
  Let $\rho:WB_n\to GL_n({\mathbb C})$ be an $n-$dimensional representation of the group $WB_n,$ $n\geqslant 3,$ such that its restriction onto $B_n$ is equivalent to the specialization of the standard representation $\tau_n(z)$ for some value $z\in {\mathbb C}^*$ with $z\neq 1$ (see Definition 2.3 and Theorem 2.5(2)).\\
  
 We will select a basis $\mathcal{B}=\{v_0, \dots ,v_{n-1}\}$ of $V,$ such that in this basis $\rho|_{B_n}=\tau_n(z).$ To find $\rho,$ it remains to determine the action of the generators $\alpha_i$ of the subgroup $S_n$ of $WB_n$ in the basis $\mathcal{B}.$\\

 To emphasize the cyclic nature of the basis $\mathcal{B},$ it is convenient to introduce an extra element $\sigma_0$ of the subgroup $B_n,$  which can be used as   a redundant generator:\\
 
$B_n=\left\langle \begin{tabular}{l|l}
&$\sigma_i \sigma_j=\sigma_j\sigma_i$ for $||i-j||\geqslant 2,$\\
$\sigma_0, \ \sigma_1,\ \sigma_2,  \dots , \sigma_{n-1}$ &$\sigma_i \sigma_{i+1}\sigma_i=\sigma_{i+1}\sigma_i \sigma_{i+1}$\\
&$\sigma_{0}=\theta\sigma_{n-1}\theta^{-1}$
\end{tabular}\right\rangle, $\\

where $\theta$ is defined in Lemma 2.2, $i, \ j \in \mathbb{Z}_n,$ and the norm of  $m$ in $\mathbb{Z}_n$ is defined by $||m||=\min\{m, n-m\}$ (here we identify the element $m$ of $\mathbb{Z}_n$ with a number $m\in\{0, \ 1, \dots, n-1\}\subseteq \mathbb{Z}$). \\

In addition to the   notations introduced in Chapter 2, we will use \\
$G_0=\rho(\sigma_0)$ and \\$C_0=\rho(\sigma_0)-I$\\

In the following lemma (Lemma 4.3) we will list some properties of the selected basis that we will use throughout sections 4 and 5. Each of the  statements in Lemma 4.3 can be easily obtained from the explicit matrices of $\tau_n(z)(\sigma_i)$ (see Definition 2.3.) The detailed discussion of the construction of this basis can be found in \cite{S}. As one can see, both images and kernels of the operators $C_i$ are conveniently written in this basis.\\

{\bf Lemma 4.3.} If $z\neq 1$ then:

(a) $v_i=Tv_{i-1}$ for $ i \in \mathbb{Z}_n,$ $i\neq 0;$ 

(b) $Tv_{n-1}=z^{n-1}v_0;$

(c) $Im(C_i)=span\{v_{i-1},v_i\}$ and \\$Ker(C_i)=span\{v_{0},\dots,v_{i-2},v_{i+1}, \dots v_{n-1}\},$ where indices are taken modulo $n.$

(d) $Im(C_{i-1})\cap Im(C_i)=span\{v_{i-1}\},$ indices are taken modulo $n.$

\vskip 0.5cm

 By Lemma 2.11(e),  it is enough to determine  $\rho(\alpha_1)$ in our selected basis $v_0,\dots ,v_{n-1}.$\\

  To classify all extensions for $n\geqslant 4,$ we  will  consider separately cases for $n\geqslant 5$ (Lemmas 4.4, 4.5, 4.6) and $n=4$ (Lemma 4.7).  The main result  for all  $n\geqslant 4$ is proven in Theorem 4.11.\\
  
{\bf Lemma 4.4.} For $n\geqslant 5$

1) Each of the vectors $v_3, v_4, \dots, v_{n-2}$ is an eigenvector of $A_1.$

2) All vectors $ v_3, v_4, \dots, v_{n-2}$  have the same  eigenvalue $k,$ where $k=1$ or $k=-1.$\\

\begin{proof}1)  For each $j=3, \dots ,n-2,$ by Lemma 4.3(d), \\ $v_j\in Im(C_{j})\cap Im(C_{j+1}).$ Thus, there exist vectors $x_j$ and $y_{j+1}$ such that $v_j=C_jx_j=C_{j+1}y_{j+1}.$

By Lemma 2.11(b), $A_1v_j=A_1C_jx_j=C_jA_1x_j\in Im(C_{j}).$ \\Similarly, $A_1v_j=A_1C_{j+1}y_{j+1}=C_{j+1}A_1y_{j+1}\in Im(C_{j+1}),$ and hence 

$A_1v_j\in Im(C_{j})\cap Im(C_{j+1})=span\{v_{j}\}$\\

2) We will show that if the set $\{ v_3, v_4, \dots, v_{n-2}\}$ consists of more than one vector ($n\geqslant 6$), then any two consecutive vectors have the same eigenvalue.\\

  Let $A_1v_j=k_jv_j$ for $j=3,4,\dots ,n-2,$ where $k_j\in  \{-1,1\}$ (See Lemma 2.11(f)).

For every $j=3, 4, \dots , n-3$ consider $A_2v_{j+1}.$ By using Lemma 2.11(e) and Lemma 4.3(a), we have\\
$A_2v_{j+1}=TA_1T^{-1}v_{j+1}=T(A_1v_j)=k_jTv_j=k_jv_{j+1}.$\\

Then, similarly to the proof of Lemma 3.6(2), \\

$k_{j+1}v_{j+1}=A_1v_{j+1}=A_2^{-1}A_1^{-1}A_2A_1A_2v_{j+1}=k_jv_{j+1},$ so $k_j=k_{j+1} ,$ since $v_{j+1}\neq 0.$

\end{proof}

{\bf Lemma 4.5.} For $n\geqslant 5$ the vector $v_2$ is an eigenvector of $A_1,$  and  its eigenvalue is equal to the eigenvalue $k$ of the vectors $v_3, v_4, \dots , v_{n-2}.$

\begin{proof} By Lemma 4.3(c),  $C_1v_3=0$ and $C_2v_3=0.$ Thus, \\
$A_2v_3=A_2(I+C_1)v_3=A_2(I+C_1)(I+C_2)v_3.$ \\By (8) and Lemma 4.4,\\
$A_2(I+C_1)(I+C_2)v_3=(I+C_1)(I+C_2)A_1v_3=k(I+C_1)(I+C_2)v_3=$\\
$=k (I+C_1)v_3=kv_3.$\\

Now, by Lemma 2.11(e) and Lemma 4.3(a),\\ $A_1v_2=T^{-1}A_2Tv_2=T^{-1}(A_2v_3)=k T^{-1}v_3=k v_2.$

\end{proof}

{\bf Lemma 4.6.} For $n\geqslant 5$ the vector $v_{n-1}$ is an eigenvector of $A_1,$ and its eigenvalue is equal to the eigenvalue $k$ of the vectors $v_2, v_3, \dots , v_{n-2}.$

\begin{proof} 

By Lemma 4.3(c), $C_1v_{n-1}=0$ and $C_2v_{n-1}=0,$  so \\$(I+C_1)^{-1}v_{n-1}=v_{n-1}$ and 
$(I+C_2)^{-1}v_{n-1}=v_{n-1}.$ Rewriting (8) as  $\alpha_i= \sigma_{i+1}^{-1}\sigma_i ^{-1}\alpha_{i+1}\sigma_i \sigma_{i+1},$ and using Lemma 2.11(e), Lemma 4.3(a), and Lemma 4.4,  we get\\

$A_1v_{n-1}=(I+C_2)^{-1}(I+C_1)^{-1}A_2(I+C_1)(I+C_2)v_{n-1}=$

$=(I+C_2)^{-1}(I+C_1)^{-1}A_2(I+C_1)v_{n-1}=(I+C_2)^{-1}(I+C_1)^{-1}A_2v_{n-1}=$

$=(I+C_2)^{-1}(I+C_1)^{-1}TA_1T^{-1}v_{n-1}=(I+C_2)^{-1}(I+C_1)^{-1}T(A_1v_{n-2})=$
$=k (I+C_2)^{-1}(I+C_1)^{-1}Tv_{n-2} =k(I+C_2)^{-1}(I+C_1)^{-1} v_{n-1}=k v_{n-1}.$

\end{proof}

   {\bf Lemma 4.7.} For $n=4$ the vectors $v_2$ and $v_3$ are the eigenvectors of $A_1$ with the same eigenvalue $k,$ where $k=1$ or $k=-1.$

\begin{proof} By Lemma 2.11(b),  $A_1C_3=C_3A_1,$  and  by Lemma 4.3(c), \\$Im(C_3)=span\{v_2,v_3\}.$ Thus,\\
$A_1v_2\in span\{v_2,v_3\}$ and  $A_1v_3\in span\{v_2,v_3\}.$\\

Suppose \\
$A_1v_2=kv_2+mv_3,$ where $k$ and $m$ are not both equal to zero. 

Then\\ $A_1v_3=A_1(I+C_3)v_2=(I+C_3)A_1v_2=(I+C_3)(kv_2+mv_3)=kv_3+mzv_2.$

Moreover, by (1), \\
$v_2=A_1^2v_2=A_1(kv_2+mv_3)=k(kv_2+mv_3)+m(kv_3+mzv_2)=(k^2+zm^2)v_2+2kmv_3,$ and hence\\
$km=0$ and $k^2+zm^2=1.$\\

We will show  that $m=0, $  and hence, $A_1v_2=kv_2$ and $A_1v_3=kv_3,$ where $k^2=1.$\\

 Suppose $m\neq 0.$ Then $k=0$ and $zm^2=1.$ Hence,\\ $A_1v_2=mv_3$ and $A_1v_3=mzv_2.$\\
 
 By Lemma 2.11(e) and Lemma 4.3 parts (b) and (a), \\$A_2v_0=TA_1T^{-1}v_0=\frac{1}{z^3}T(A_1v_3)=\frac{1}{z^3}T(mzv_2)=\frac{m}{z^2}v_3$\\

By group identity (7), together with $zm^2=1,$ we get \\
$A_1A_2v_1=A_1A_2(I+C_1)v_0=(I+C_2)A_1A_2v_0=\frac{m}{z^2}(I+C_2)A_1v_3=$\\
$=\frac{m\cdot mz}{z^2}(I+C_2)v_2=\frac{1}{z^2}\cdot zv_1=\frac{1}{z}v_1.$\\

By Lemma 2.2(d),  $(A_1A_2)^3=I,$ so $z^3=1.$ Thus, since $z\neq 1,$ the case $m\neq 0$ only possible when $z=\omega$ or $z=\omega^2,$ where $\omega$ is the primitive cubic root of unity. \\

Suppose now that $z=\omega$ or $z=\omega^2.$ By Lemma 2.11(c), $A_1A_2C_1=C_2A_1A_2$ and and by Lemma 4.3(c),\\$C_1v_3=0.$ Thus, by using Lemma 2.11(e) and Lemma 4.3 parts (a) and (b), we have\\

$0=A_1A_2C_1v_3=C_2A_1A_2v_3=C_2A_1(TA_1T^{-1}v_3)=C_2A_1T(A_1v_2)=$

$=C_2A_1T(mv_3)=mC_2A_1(z^3v_0)=mC_2(A_1v_0)$\\

Since $m\neq 0,$ we get that $A_1v_0\in Ker(C_2)=span\{v_0,v_3\}$\\

Suppose $A_1v_0=av_0+bv_3.$  Then from\\

$v_0=A_1^2v_0=A_1(av_0+bv_3)=a(av_0+bv_3)+b(zmv_2)$ it follows that $b=0$ and $a^2=1.$ Hence,\\
$A_1v_0=av_0.$\\

From $A_1A_2v_1=\frac{1}{z}v_1,$ and Lemma 2.11(e) and Lemma 4.3(a), we have
 \\$A_1v_1=zA_2v_1=zT(A_1v_0)=zT(av_0)=zav_1,$ and \\

$v_1=A_1^2v_1=z^2a^2v_1=z^2v_1.$ That gives $z^2=1,$ a contradiction with $z=\omega$ or $z=\omega^2.$\\

\end{proof}

The following corollary summarizes the results of Lemmas 4.4, 4.5, 4.6 and 4.7.\\

{\bf Corollary 4.8.} Let $\rho:WB_n\to GL_n({\mathbb C})$ be an $n-$dimensional representation of $WB_n,$ $n\geqslant 4,$ such that its restriction $\rho|_{B_n}$ onto $B_n$ is equivalent to the specialization of the standard representation $\tau_n(z)$ for some value $z\in {\mathbb C}^*,$  $z\neq 1.$  Then $\rho(\alpha_1)v_j=kv_j \ \  \forall j=2, \dots, n-1,$ where $k\in \{-1,1\}.$\\

Similarly to Chapter 3, we will consider only those representations $\rho,$ for which  $\rho(\alpha_1)v_j=v_j, \ j=2, 3, \dots, n-1.$  The  classification theorem for $n-$dimensional extensions of the standard representations will follow then from  this specific case.\\

{\bf Lemma 4.9.} Suppose $\rho:WB_n\to GL_{n}({\mathbb C})$ for  $n\geqslant 4,$ such that  $\rho|_{B_n}$ is equivalent to $\tau_n(z),$ where $z\in \mathbb{C^*},$  $z\neq 1.$ Suppose that  $\rho(\alpha_1)v_j=v_j$ for $j=2, 3, \dots , n-1.$ \\

Then $\exists \lambda\in \mathbb{C}^*,$ such that $A_1v_0=\lambda v_1$ and $A_1v_1=\frac{1}{\lambda}v_0.$

\begin{proof} 

By Lemma  4.3(c), $Im(C_1)=span\{v_0,v_1\}$ and \\$Im(C_2)=span\{v_1,v_2\},$ and by Lemma 2.11(c) $A_1A_{2}C_1=C_{2}A_1A_{2},$ so both\\

$A_1A_{2}v_0\in span\{v_1,v_2\}$ and $A_1A_{2}v_1\in span\{v_1,v_2\}$\\

By Lemma 2.11(e) and Lemma 4.3 (b), \\$A_2v_0=TA_1T^{-1}v_0=\frac{1}{z^{n-1}}TA_1v_{n-1}=\frac{1}{z^{n-1}}Tv_{n-1}=v_0,$ so\\
$A_1v_0=A_1A_{2}v_0\in span\{v_1,v_2\}.$ \\

Let $A_1v_0=\lambda v_1+\mu v_2.$\\

Again, by Lemma 2.11(e) and Lemma 4.3 (a), \\$A_2v_1=TA_1T^{-1}v_1=T(A_1v_0)=T(\lambda v_1+\mu v_2)=\lambda v_2+\mu v_3,$ \\
so $A_1A_{2}v_1=A_1(\lambda v_2+\mu v_3)=\lambda v_2+\mu v_3\in span\{v_1,v_2\},$ thus \\
$\mu=0.$ 

Hence, $A_1v_0=\lambda v_1$ with  $\lambda\neq 0$ since $A_1$ is invertible. Additionally, from $A_1^2=I,$ it follows that $A_1v_1=\frac{1}{\lambda}v_0.$\\

\end{proof}

{\bf Lemma 4.10.} Suppose $\rho:WB_n\to GL_{n}({\mathbb C})$ for  $n\geqslant 4,$ such that  $\rho|_{B_n}$ is equivalent to $\tau_n(z),$ where $z\in \mathbb{C^*},$  $z\neq 1.$ Suppose that  $\rho(\alpha_1)v_j=v_j$ for $j=2, 3, \dots , n-1.$

Then $\rho$ is equivalent to $\widetilde{\tau}_n (z,\lambda)$ for some $ \lambda\in \mathbb{C}^*.$\\

\begin{proof} Select a basis $\mathcal{B}=\{v_0, \dots ,v_{n-1}\}$ of $V,$ such that in this basis $\rho(\sigma_i)=\tau_n(z)(\sigma_i),$ $i=1, \dots n-1,$  where $z\neq 1, \  z\in \mathbb{C^*}.$ Then by Lemma 4.9, there exists $ \lambda\in \mathbb{C}^*,$ such that 
$\rho(\alpha_1)=\widetilde{\tau}_n(z,\lambda)(\alpha_1).$	\\

By Lemma 4.3, parts (a) and (b), the matrix of $T$ in  basis $\mathcal{B} $ is \\

$T=\left( \begin{array}{ccccc}
	
	0&0&\dots&0&z^{n-1}\\
	1&0&\dots&0&0\\ 
	0&1&\dots&0&0\\
	\vdots&\vdots&\ddots&\vdots&\vdots \\
	0&0&\dots&1&0
\end{array} \right)$\\ 

(1s are located below the main diagonal).\\

Now,  from Lemma 2.11(e), the definition of $\widetilde{\tau}_n (z,\lambda)$ (Definition 4.1) and the explicit expressions for the matrices $T$ and $T^{-1},$ we have 
 $\rho(\alpha_i)=	\widetilde{\tau}_n (z,\lambda)(\alpha_i)$ for $i=2, \dots , n-1,$  hence, $\rho$ is equivalent to $\widetilde{\tau}_n (z,\lambda).$\\
	
\end{proof}

{\bf Theorem 4.11.} Let $n\geqslant 4.$ Let $\rho:WB_n\to GL_n(\mathbb{C})$ be an $n-$dimensional representation of the welded braid group $WB_n,$ such that its restriction $\rho|_{B_n}$ onto the braid group $B_n$ is equivalent to a specialization of a standard representation $\tau_n(z)$ for some $z\neq 1,\ z\in\mathbb{C}^*.$ 

Then there exists a unique pair $(\lambda, k)$ where $\lambda \in\mathbb{C}^*$ and $k\in \{-1,1\}$ such that $\rho$ is equivalent to $X_n(1,k)\otimes\widetilde{\tau}_n(z,\lambda).$

\begin{proof}   Select a basis $\mathcal{B}=\{v_0, \dots ,v_{n-1}\}$ of $V,$ such that in this basis $\rho(\sigma_i)=\tau_n(z)(\sigma_i),$ $i=1, \dots n-1,$  where $z\neq 1.$\\

By Corollary 4.8,   there exists $k\in\{-1,1\},$ such that \\$\rho(\alpha_1)v_j=kv_j \ \  \forall j=2, \dots, n-1.$ \\

Consider $\hat{\rho}=(X_n(1,k))^{-1}\otimes {\rho},$ where $X_n(1,k)$ is a one-dimensional representation of $WB_n.$ Then $\hat{\rho}$ is an $n-$dimensional representation of $WB_n,$ such that $\hat{\rho}_{|B_n}$ is equivalent to $\tau_n(z)$ and \\$\hat{\rho}(\alpha_1)v_j=v_j$ for all $j=2,3,\dots, n-1.$\\

By Lemma 4.10, $\hat{\rho}$ is equivalent to $\widetilde{\tau}_n (z,\lambda)$ for some $ \lambda\in \mathbb{C}^*,$
thus, $\rho=X_n(1,k)\otimes \hat{\rho}$ is equivalent to $X_n(1,k)\otimes\widetilde{\tau}_n (z,\lambda).$\\

 Uniqueness of the pair $(\lambda, k)$ immediately follows from the fact that  \\ $X_n(1,k_1)\otimes\widetilde{\tau}_n(z,\lambda_1)\neq X_n(1,k_2)\otimes\widetilde{\tau}_n(z,\lambda_2)$ whenever $(\lambda_1, k_1)\neq(\lambda_2, k_2).$\\Since $z\neq 1,$ the restriction $\rho|_{B_n}$ is irreducible, so, by Lemma 2.9, $X_n(1,k_1)\otimes\widetilde{\tau}_n(z,\lambda_1)$ is not equivalent to $X_n(1,k_2)\otimes\widetilde{\tau}_n(z,\lambda_2).$ 

\end{proof}

\section{Extensions of the standard irreducible representation for $n=3.$}

In this section we will classify all extensions of an irreducible standard representation for $n=3.$ We will continue to use the notations introduced in section 4. The main classification theorem (Theorem 5.9) for $n=3$ is proven at the end of this section. \\

{\bf Definition 5.1.} For $z,\lambda\in \mathbb{C}^*,$  define 

$\psi_{3}(z,\lambda;x):WB_3\to GL_3(\mathbb{C})$ by\\

$$\psi_{3}(z,\lambda;x) (\sigma_1)=\left( \begin{array}{rrr}
0&z&0\\
 
1&0&0\\
0&0&1
\end{array} \right),
\psi_{3}(z,\lambda;x) (\sigma_2)=\left( \begin{array}{rrr}
1&0&0\\
 
0&0&z\\
0&1&0
\end{array} \right),$$ \\

$$\psi_{3}(z,\lambda;x) (\alpha_1)=\left( \begin{array}{ccc}
a&\frac{b}{x}&\frac{c}{x^2}\\
 &&\\
xb&c&\frac{a}{x}\\
&&\\
x^2c&xa&b
\end{array} \right),
\psi_{3}(z,\lambda;x) (\alpha_2)=\left( \begin{array}{ccc}
b&\frac{c}{x}&\frac{a}{x^2}\\
&&\\
 
xc&a&\frac{b}{x}\\
&&\\
x^2a&xb&c
\end{array} \right),$$ \\

where $x^3=\frac{1}{z^2}$ and \\ 

$\left\{ \begin{array}{lcl}
a=a(\lambda)&=&\frac{\omega}{3}\lambda+\frac{1}{3}+\frac{\omega^2}{3}\lambda^{-1}\\
&&\\
b=b(\lambda)&=&\frac{1}{3}\lambda+\frac{1}{3}+\frac{1}{3}\lambda^{-1}\\
 &&\\ 
c=c(\lambda)&=&\frac{\omega^2}{3}\lambda+\frac{1}{3}+\frac{\omega}{3}\lambda^{-1}
\end{array} \right.$\\

\vskip 0.5cm

{\bf Remark 5.2.} 1) The role of the parameter $x$ in the notation of $\psi_{3}(z,\lambda;x)$ is different from the parameters $z$ and $\lambda$: for every fixed $z\in \mathbb{C}^*$ there are three choices for $x,$ each equal to one of the cubic roots of $ \frac{1}{z^2}, $ and each fixed pair $(z,x)$ defines a one-parameter family of 3-dimensional representations of $WB_3,$ parametrized by $\lambda \in\mathbb{C}^*;$\\
2) The correspondence $\lambda \mapsto (a(\lambda),b(\lambda),c(\lambda))$ is injective; \\
3) The point $(0,1,0)$ corresponds to the value of $\lambda=1.$\\

{\bf Lemma 5.3.}   

1) For all $z,\lambda\in \mathbb{C}^*,$ the map $\psi_{3}(z,\lambda;x):WB_3\to GL_3(\mathbb{C}) $ is a representation of $WB_3,$ which is an extension of the specialization of the  standard representation $\tau_3(z)$  of the braid group  $B_3.$ If  $z\neq 1,$ then $\psi_{3}(z,\lambda;x)$ is irreducible.\\

2) For all $z_1,z_2,\lambda_1,\lambda_2 \in \mathbb{C}^*:$\\
a) if $z_1\neq z_2,$ then $\psi_{3}(z_1,\lambda_1;x_1)$ and $\psi_{3}(z_2,\lambda_2;x_2)$ are  not equivalent;\\
b) if $z_1=z_2=z\neq 1$ and $(\lambda_1,x_1)\neq (\lambda_2,x_2)$ then $\psi_{3}(z,\lambda_1;x_1)$ and $\psi_{3}(z,\lambda_2;x_2)$ are  not equivalent;\\

3) For all $z_1,z_2,\lambda_1,\lambda_2 \in \mathbb{C}^*:$\\
a) if $z_1\neq z_2,$ then $\psi_{3}(z_1,\lambda_1;x_1)$ and $\widetilde{\tau}_3 (z_2,\lambda_2)$ are  not equivalent;\\
b) if $z_1=z_2=z\neq 1$ and $\lambda_1\neq 1,$ then $\psi_{3}(z,\lambda_1;x_1)$ and $\widetilde{\tau}_3 (z,\lambda_2)$ are  not equivalent;  for $\lambda_1=1$ the representation $\psi_{3}(z,1;x)$ is equal to $\widetilde{\tau}_3 (z,x).$\\

\begin{proof}  It is easy to check that all defining relations (1)--(8) are satisfied, and $\psi_{3}(z,\lambda;x)(\sigma_i)=\tau_n(z)(\sigma_i)$ for \ \ $ i=1,2.$ By Theorem 2.5, part 2), $\tau_n(z)$ is irreducible for $z\neq 1,$ so $\psi_{3}(z,\lambda;x)$ is irreducible, which proves part 1).\\
	
Like in  Lemma 4.2, part 3),  the restrictions $\tau_3(z)$ of the representations $\psi_{3}(z,\lambda;x)$ and $\widetilde{\tau}_3 (z,\lambda)$ onto $B_3$ are non-equivalent for distinct values of $z,$ from where parts 2a) and 3a) immediately follow. 

For $z_1=z_2=z\neq 1,$ the restriction  $\tau_3(z)$ is irreducible, so  by Lemma 2.9,  it is enough to see that the corresponding  representations are not equal. Indeed, for part 2b), by Remark 5.2, part 2),$\psi_{3}(z,\lambda_1;x_1)(\alpha_1)\neq \psi_{3}(z_2,\lambda_2;x_2) (\alpha_1)$ for all $\lambda_1\neq \lambda_2,$ and $\psi_{3}(z,\lambda;x_1)(\alpha_1)\neq \psi_{3}(z_2,\lambda;x_2) (\alpha_1)$ for $x_1\neq x_2.$ And for part 3b), if $\lambda_1\neq 1,$ then by Remark 5.2, parts 2) and 3), $(a,b,c)\neq (0,1,0),$ so $\psi_{3}(z,\lambda_1;x_1)\neq \widetilde{\tau}_3 (z,\lambda_2).$ The last statement of part 3b) follows from the Definitions 4.1 and 5.1.

\end{proof}

 Similarly to Chapter 4,  let $\rho:WB_3\to GL_3({\mathbb C})$ be a $3-$dimensional representation of the group $WB_3,$  such that its restriction onto $B_3$ is equivalent to the specialization of the standard representation $\tau_3(z)$ for some value $z\in {\mathbb C}^*$ with $z\neq 1.$\\

Select a basis $\mathcal{B}=\{v_0,v_1,v_2\}$ of $V,$ such that in this basis \\$\rho|_{B_3}=\tau_3(z).$ We will determine the action of the generators $\alpha_1$ and $\alpha_2$  in the basis $\mathcal{B}.$\\

{\bf Lemma 5.4.} Suppose $n=3.$ Consider $A=A_1A_2=\rho(\alpha_1)\rho(\alpha_2).$\\ Then the matrix of $A$ written in the basis $\mathcal{B}=\{v_0,v_1,v_2\}$ is\\

$A=\left( \begin{array}{ccc}

0&0&\frac{1}{x^2}\\
 
x&0&0\\
0&x&0
\end{array} \right),$ where $x\in \mathbb{C}^*.$

\begin{proof}
Similarly to the proof of Lemma 4.9, $Im(C_1)=span\{v_0,v_1\}$ and $Im(C_2)=span\{v_1,v_2\}$ by Lemma  4.3(c),  and $A_1A_{2}C_1=C_{2}A_1A_{2}$  by Lemma 2.11(c). Hence $A_1A_{2}v_0\in span\{v_1,v_2\}.$ \\
 
 Let $Av_0=A_1A_{2}v_0=x_1v_1+x_2v_2.$ Then, using group identity (7), \\
 $Av_1=A_1A_{2}v_1=A_1A_{2}(I+C_1)v_0=(I+C_2)A_1A_{2}v_0=$\\
 $=(I+C_2)(x_1v_1+x_2v_2)=x_1v_2+x_2zv_1$
 
 By Lemma 4.3(c), $v_2\in Ker(C_1),$ so by Lemma 2.11(c),\\
 $0=A_1A_2C_1v_2=C_2A_1A_2v_2,$ so $Av_2=A_1A_2v_2\in Ker(C_2)=span\{v_0\}.$  Hence  $Av_2=x_0v_0$ for some $x_0\in\mathbb{C}^*.$
 \\
Writing the matrix of   $A$   in the basis $\{v_0,v_1,v_2\},$  we obtain \\ 

$A=\left( \begin{array}{ccc}

0&0&x_0\\
 
x_1&x_2z&0\\
x_2&x_1&0
\end{array} \right)$\\

By Lemma 2.2(d),  $A^3=I,$ and by direct calculation we have\\

$A^3=\left( \begin{array}{ccc}

x_0x_1^2&*&x_0^2x_2\\
 
*&*&*\\
*&*&*
\end{array} \right),$\\

from which we have  $x_0x_1=1$ and $x_0^2x_2=0.$ Hence, $x_0\neq 0, \ x_1\neq 0, $ and so $x_2=0.$ Denoting $x_1=x, $ we get $x_0=\frac{1}{x^2},$ and  the required result immediately follows.\\

\end{proof}

 Unlike  the general case, for $n=3$ the vector $v_2$ may or may not be the eigenvector of $\rho(\alpha_1).$ We will consider these two cases separately. \\
 
 For the case when $v_2$ is an  eigenvector of $\rho(\alpha_1),$ similarly to the general case, by tensoring with the appropriate one-dimensional representation, it is enough to consider only representations with the eigenvalue 1. In the next lemma we prove the result similar to the result of   Lemma 4.9 for $n=3.$\\
 
{\bf Lemma 5.5.} Suppose $\rho:WB_3\to GL_{3}({\mathbb C}),$  such that  $\rho|_{B_3}$ is equivalent to $\tau_3(z),$ where $z\in \mathbb{C^*},$  $z\neq 1.$ Suppose that  $\rho(\alpha_1)v_2=v_2.$\\

Then $\exists \lambda\in \mathbb{C}^*,$ such that $\rho(\alpha_1)v_0=\lambda v_1$ and $\rho(\alpha_1)v_1=\frac{1}{\lambda}v_0.$

\begin{proof} By Lemma 5.4, $A_1A_{2}v_1=xv_2$ for some $x\in\mathbb{C}^*.$ Since $A_1^2=I,$ we have \\$A_{2}v_1=xA_{1}v_2=xv_2.$ \\

Let $\lambda=x\in \mathbb{C}^*.$ Then, by Lemma 2.11(e) and Lemma 4.3(a),  \\
$A_1v_0=T^{-1}(A_2v_1)=T^{-1}(\lambda v_2)=\lambda v_1$ and again, since $A_1^2=I,$ \\
$A_1v_1=\frac{1}{\lambda}v_0.$\\
\end{proof}

Next, we will consider the case when $v_2$ is not an  eigenvector of $\rho(\alpha_1).$\\

{\bf Lemma 5.6.} Suppose $\rho:WB_3\to GL_{3}({\mathbb C}),$  such that  $\rho|_{B_3}$ is equivalent to $\tau_3(z),$ where $z\in \mathbb{C^*},$  $z\neq 1.$ Suppose that the vector $v_2$ is {\bf not} an eigenvector of $\rho(\alpha_1).$\\ Then the matrices for $\rho(\alpha_1)$ and $\rho(\alpha_2)$ in the basis  $\{v_0,v_1,v_2\}$ are

$$\rho(\alpha_1)=k\cdot \left( \begin{array}{ccc}
a&\frac{b}{x}&\frac{c}{x^2}\\
 &&\\
xb&c&\frac{a}{x}\\
&&\\
x^2c&xa&b
\end{array} \right),\hskip 0.5cm
\rho(\alpha_2)=k\cdot \left( \begin{array}{ccc}
b&\frac{c}{x}&\frac{a}{x^2}\\
&&\\
 
xc&a&\frac{b}{x}\\
&&\\
x^2a&xb&c
\end{array} \right),$$ \\

where $k\in\{-1,1\},$ $x^3=\frac{1}{z^2},$  $(a,b,c)$ satisfy \\ 

$\left\{ \begin{array}{lcr}
a^2+b^2+c^2&=&1\\

a+b+c&=&1

\end{array} \right. ,$\\

and $(a,b,c)\neq (0,1,0).$\\

\begin{proof} By  Lemma 5.4, \\

$A_1A_2v_0=xv_1,$

$A_1A_2v_1=xv_2,$

$A_1A_2v_2=\frac{1}{x^2}v_0$, so\\

$A_2v_0=xA_1v_1,$

$A_2v_1=xA_1v_2,$

$A_2v_2=\frac{1}{x^2}A_1v_0.$\\

Let $A_1v_0=av_0+bv_1+cv_2.$ \\

Then $A_2v_2=\frac{a}{x^2}v_0+\frac{b}{x^2}v_1+\frac{c}{x^2}v_2,$ hence, by Lemma 2.11(e) and Lemma 4.3 parts (a) and (b),\\

 $A_1v_1=T^{-1}(A_2v_2)=T^{-1}(\frac{a}{x^2}v_0+\frac{b}{x^2}v_1+\frac{c}{x^2}v_2)=\frac{a}{z^2x^2}v_2+\frac{b}{x^2}v_0+\frac{c}{x^2}v_1.$\\
 
 Similarly, \\
 
 $A_1v_2=\frac{1}{x}A_2v_1=\frac{1}{x}T(A_1v_0)=\frac{1}{x}T(av_0+bv_1+cv_2)=\frac{a}{x}v_1+\frac{b}{x}v_2+\frac{z^2c}{x}v_0.$\\
 
 Writing down the matrices for $A_1$ and $A_2,$ we obtain:\\
 
$A_1= \left( \begin{array}{ccc}
a&\frac{b}{x^2}&\frac{z^2c}{x}\\
&&\\
 
b&\frac{c}{x^2}&\frac{a}{x}\\
&&\\
c&\frac{a}{x^2z^2}&\frac{b}{x}
\end{array} \right)$  and $A_2=TA_1T^{-1}= \left( \begin{array}{ccc}
\frac{b}{x}&z^2c&\frac{a}{x^2}\\
&&\\
 
\frac{c}{x}&a&\frac{b}{x^2}\\
&&\\
\frac{a}{xz^2}&b&\frac{c}{x^2}
\end{array} \right),$\\
\vskip 0.3cm

 where $T=\left( \begin{array}{rrr}
0&0&z^2\\
 
1&0&0\\
0&1&0
\end{array} \right).$\\

By direct calculations we get\\

$(A_1A_2C_1-C_2A_1A_2)v_0=\left[ \begin{array}{c}
*\\
*\\
\frac{c^2(x^3z^2-1)}{x^3}
\end{array} \right]$  and \\

$(C_1A_2A_1-A_2A_1C_2)v_2=\left[ \begin{array}{c}
\frac{a^2(x^3z^2-1)}{zx^4}\\
*\\
*
\end{array} \right].$\\

Since $v_2$ is not an eigenvector of $A_1,$ either $a\neq 0$ or $c\neq 0$ (or both). Thus, 
 $ x^3z^2-1=0,$ and so $x^3=\frac{1}{z^2}.$\\
 
We will substitute $\frac{1}{x^3}$ for $z^2$  in $A_1$ and $A_2$  and, for the sake of symmetry, replace
 $b$ by $xb$ and  $c$ by $x^2c.$ We obtain\\
 
 $A_1=\left( \begin{array}{ccc}
a&\frac{b}{x}&\frac{c}{x^2}\\
&&\\
 
xb&c&\frac{a}{x}\\
&&\\
x^2c&xa&b
\end{array} \right),\hskip 1cm
A_2=\left( \begin{array}{ccc}
b&\frac{c}{x}&\frac{a}{x^2}\\
&&\\
 
xc&a&\frac{b}{x}\\
&&\\
x^2a&xb&c
\end{array} \right).$ \\

By direct calculation we have\\

$I=A_1^2=\left( \begin{array}{ccc}
a^2+b^2+c^2&\frac{1}{x}(ab+bc+ac)&\frac{1}{x^2}(ab+bc+ac)\\
&&\\
 
x(ab+bc+ac)&a^2+b^2+c^2&\frac{1}{x}(ab+bc+ac)\\
&&\\
x^2(ab+bc+ac)&x(ab+bc+ac)&a^2+b^2+c^2
\end{array} \right),$\\

which gives\\

 $\left\{ \begin{array}{lcr}
a^2+b^2+c^2&=&1\\

ab+bc+ac&=&0

\end{array} \right. $\\

This system of equations is  equivalent to \\

$\left\{ \begin{array}{lcr}
a^2+b^2+c^2&=&1\\

(a+b+c)^2&=&1

\end{array} \right. $ \\ 

We will rewrite this system as \\

$\left\{ \begin{array}{lcr}
a^2+b^2+c^2&=&1\\

a+b+c&=&k

\end{array} \right. ,$ where $k^2=1,$ or, alternatively, as\\

\vskip 0.3cm

$\left\{ \begin{array}{lcr}
\left(\frac{a}{k}\right)^2+\left(\frac{b}{k}\right)^2+\left(\frac{c}{k}\right)^2&=&1\\
&&\\

\left(\frac{a}{k}\right)+\left(\frac{b}{k}\right)+\left(\frac{c}{k}\right)&=&1

\end{array} \right. ,$ where $k\in\{-1,1\}.$\\

\vskip 0.3cm

Replacing $\left(\frac{a}{k}\right),\ \left(\frac{b}{k}\right)$ and $\left(\frac{c}{k}\right)$ by $a,\ b$ and $c $ respectively, we obtain that \\

$$A_1=k\cdot \left( \begin{array}{ccc}
a&\frac{b}{x}&\frac{c}{x^2}\\
 &&\\
xb&c&\frac{a}{x}\\
&&\\
x^2c&xa&b
\end{array} \right),\hskip 1cm
A_2=k\cdot \left( \begin{array}{ccc}
b&\frac{c}{x}&\frac{a}{x^2}\\
&&\\
 
xc&a&\frac{b}{x}\\
&&\\
x^2a&xb&c
\end{array} \right),$$ \\

where $k\in\{-1,1\}$  and  \\ 

$\left\{ \begin{array}{lcr}
a^2+b^2+c^2&=&1\\

a+b+c&=&1

\end{array} \right. .$\\

Since $v_2$ is not an eigenvector of $A_1,$ at least one of $a$, $c$ is non-zero, hence $(a,b,c)\neq (0,1,0).$\\

\end{proof}

{\bf Lemma 5.7.} The complex curve $\left\{ \begin{array}{lcr}
{\bf x}^2+{\bf y}^2+{\bf z}^2&=&1\\

{\bf x}+{\bf y}+{\bf z}&=&1

\end{array} \right. $ \\can be parametrized by \\

$\left\{ \begin{array}{lcl}
{\bf x}={\bf x}(t)&=&\frac{\omega}{3}t+\frac{1}{3}+\frac{\omega^2}{3}t^{-1}\\
&&\\
{\bf y}={\bf y}(t)&=&\frac{1}{3}t+\frac{1}{3}+\frac{1}{3}t^{-1}\\
 &&\\ 
{\bf z}={\bf z}(t)&=&\frac{\omega^2}{3}t+\frac{1}{3}+\frac{\omega}{3}t^{-1}
\end{array} \right.$ ,$\hskip 1cm  t\in \mathbb{C}^*$\\

(Here $\omega$ denotes the primitive cubic root of unity).

\begin{proof}

We can parametrize this curve by \\

$\left\{ \begin{array}{lcr}
{\bf x}&=&\frac{s+1}{s^2+s+1}\\
{\bf y}&=&\frac{s^2+s}{s^2+s+1}\\
  
{\bf z}&=&-\frac{s}{s^2+s+1}
\end{array} \right. , \hskip 0.5cm \begin{array}{l}{\textrm{where }}s\neq \omega, \omega ^2,\\
s\in \mathbb{C}\cup \{\infty\}\end{array}$ \\

By changing the parameter to $t=\frac{s-\omega}{s-\omega^2}, \ t\in \mathbb{C}^*,$ we obtain the required parametrization.\\

\end{proof}

By using Definition 5.1 together with Lemmas 5.6 and 5.7 and Remark 5.2 parts 2) and 3), we obtain the following result.\\

{\bf Corollary 5.8.} Suppose $\rho:WB_3\to GL_{3}({\mathbb C}),$  such that  $\rho|_{B_3}$ is equivalent to $\tau_3(z),$ where $z\in \mathbb{C^*},$  $z\neq 1.$ Suppose that the vector $v_2$ is {\bf not} an eigenvector of $\rho(\alpha_1).$\\ 

Then there exist $k\in \{-1,1\},$ $\lambda \in\mathbb{C}^*, \lambda\neq 1,$ and $x \in\mathbb{C}^*, \ x^3=\frac{1}{z^2},$ such that $\rho$ is equivalent to $X_3(1,k)\otimes\psi_3(z,\lambda;x).$\\

{\bf Theorem 5.9.}  Let $\rho:WB_3\to GL_3(\mathbb{C})$ be a $3-$dimensional representation of the welded braid group $WB_3,$  such that its restriction $\rho|_{B_3}$ onto the braid group $B_3$ is equivalent to a specialization of a standard representation $\tau_3(z)$ for some $z\neq 1,\ z\in\mathbb{C}^*.$ Then exactly one of the following is true:\\

\noindent
1)  There exists a unique pair $(\lambda, k)$ where $\lambda \in\mathbb{C}^*$ and $k\in \{-1,1\}$ such that $\rho$ is equivalent to $X_3(1,k)\otimes\widetilde{\tau}_3(z,\lambda);$\\

\noindent
2) There exists a unique triple $(\lambda, k,x)$ where $k\in \{-1,1\},$ $\lambda \in\mathbb{C}^*,$ $\lambda\neq 1,$ and $x \in\mathbb{C}^*, \ x^3=\frac{1}{z^2},$ such that $\rho$ is equivalent to \\ $X_3(1,k)\otimes\psi_3(z,\lambda;x).$\\

\begin{proof}   Select a basis $\mathcal{B}=\{v_0, v_1 ,v_{2}\}$ of $V,$ such that in this basis $\rho(\sigma_i)=\tau_3(z)(\sigma_i),$ $i=1,2.$ Then the vector $v_{2}$ either is  an eigenvector  $\rho(\alpha_1)$ (Case (1)), or $v_{2}$ is not an eigenvector  $\rho(\alpha_1)$ (Case(2)).     \\

1) The vector $v_2$ is an eigenvector of $\rho(\alpha_1).$  Let $\rho(\alpha_1)v_2=kv_2 .$ By Lemma 2.11(f), $k\in\{-1,1\}.$\\

Consider $\hat{\rho}=(X_3(1,k))^{-1}\otimes {\rho},$ where $X_3(1,k)$ is a one-dimensional representation of $WB_3.$ Then $\hat{\rho}$ is an $3-$dimensional representation of $WB_3,$ such that $\hat{\rho}_{|B_3}$ is equivalent to $\tau_3(z)$ and $\hat{\rho}(\alpha_1)v_2=v_2.$\\ By Lemma 5.5, there exists $ \lambda\in \mathbb{C}^*,$ such that 
$\hat{\rho}(\alpha_1)=\widetilde{\tau}_3(z,\lambda)(\alpha_1).$ Then, similarly to Lemma 4.10, using Lemma 2.11(e), the explicit expressions for $T$  and $T^{-1}$ in basis $\mathcal{B},$ obtained from Lemma 4.3 parts (a) and (b),  and Definition 4.1, we have $\hat{\rho}(\alpha_2)=	\widetilde{\tau}_3 (z,\lambda)(\alpha_2).$ Thus, $\hat{\rho}$ is equivalent to $\widetilde{\tau}_3 (z,\lambda),$ and hence, $\rho=X_3(1,k)\otimes \hat{\rho}$ is equivalent to $X_3(1,k)\otimes\widetilde{\tau}_3 (z,\lambda).$\\

2) The vector  $v_2$ is not an eigenvector of $\rho(\alpha_1).$ Then by Corollary 5.8, there exist $k\in \{-1,1\},$ $\lambda \in\mathbb{C}^*,$ $\lambda\neq 1,$ and $x \in\mathbb{C}^*, \ x^3=\frac{1}{z^2},$ such that $\rho$ is equivalent to $X_3(1,k)\otimes\psi_3(z,\lambda;x).$\\

Now we will show that  for $z\neq 1,$ every two distinct  extensions  $X_3(1,k)\otimes \tau',$ where $\tau'$ is either $\widetilde{\tau}_3 (z,\lambda),$  $\lambda \in\mathbb{C}^*,$ or $\psi_3(z,\lambda;x),$ $\lambda \in\mathbb{C}^*,$ $\lambda\neq 1,$ are not equivalent. 

First, since $Trace\left( (X_3(1,k)\otimes \tau')(\alpha_1)\right) =k, $ distinct values of $k$ give non-equivalent representations. Next, by Lemma 4.2 part 3b), the representations $\widetilde{\tau}_3(z,\lambda)$ are not equivalent for distinct values of $\lambda.$ By Lemma 5.3 part 2b), the representations $\psi_3(z,\lambda;x)$ are not equivalent for distinct pairs $(\lambda, x).$ And by Lemma 5.3 part 3b), the representation $\psi_3(z,\lambda;x)$ is not equivalent to $\widetilde{\tau}_3(z,\lambda')$ for $\lambda\neq 1.$\\

\end{proof}

\section{Extensions of the irreducible representations of $B_n$ of dimension less or equal to $n$}

The goal of this section is to find (up to equivalence) all extensions of every irreducible complex representation of $B_n$ of dimension $\leqslant n$ to the welded braid group    $WB_n$ for $n$ large enough.\\

All irreducible representations of the braid group $B_n$ of dimension up  to $n$ were classified  in a series of papers \cite{Formanek},\cite{Lee},\cite{S}, and \cite{FLSV} for all $n\geqslant 3.$ It has been proven in \cite{Formanek}, Theorem 23, that for $n\geqslant 7$  the only irreducible representations of $B_n$ of dimension $n-1$ are those of a Burau type, while there is  a number of exceptions for   $n\leqslant 6.$ Likewise, 
it has been proven in \cite{S}, Theorem 6.1 that for $n\geqslant 9$  the only irreducible representations of $B_n$ of dimension $n$ are those of a standard type. Additionally, it has been proven in \cite{FLSV}, Theorem 6.1 that  for $n=7$  there are no exceptional cases, while for $n\leqslant 6$ and $n=8$ there is a number of exceptions. We will formulate and prove our classification theorems (Theorem 6.2 and Theorem 6.3) only for those values of $n$ where there are no exceptional cases. The classification theorems for the other values of $n$ require the classification of the extensions of each of the exceptional cases, and is outside of the scope of this paper.\\

In the first theorem (Theorem 6.1) we prove that the welded braid group $WB_n$ has no irreducible representations of dimension $r,$ where $2\leqslant r \leqslant n-2 $ for $n\geqslant 5.$ This result is similar to the Burnside theorem for symmetric groups and is a direct consequence of it. It follows immediately from this result that for $n\geqslant 5$ any irreducible representation of $B_n$ of dimension $n-2$ can not be extended to $WB_n.$ (For the classification of $(n-2)-$dimensional representations of $B_n$ see \cite{Formanek}, Theorem 23(b).)\\

{\bf Theorem 6.1.} Suppose $n\geqslant 5.$ Let $\rho:WB_n\to GL_{r}({\mathbb C})$ be an $r-$dimensional representation of the welded braid  group $WB_n,$ where $2\leqslant r \leqslant n-2.$ Then $\rho$ is {\it reducible}.

\begin{proof} The idea of the proof of this theorem is very similar to the proof of Lemma 3.9. Consider the restriction $\rho|_{S_n}$ onto symmetric group $S_n.$ By the theorem of Burnside (\cite{Burnside}, note C, p. 468), for $n\geqslant 5,$ the symmetric group $S_n$ has no irreducible representations of dimension $r$ for $2\leqslant r\leqslant n-2,$ hence, by Maschke's theorem,  $\rho|_{S_n}$ is a direct sum of one-dimensional representations. For any one-dimensional representation $\eta$ of $S_n,$ we have $\eta(\alpha_i)=\eta(\alpha_j)$ for all $1\leqslant i,j\leqslant n-1$. Hence, $\rho(\alpha_i)=\rho(\alpha_j)$ for all $1\leqslant i,j\leqslant n-1,$ and from the group identities (7) and (1), it follows that $\rho(\sigma_i)=\rho(\sigma_{i+1})$,  for $i=1, 2,\dots,n-2.$ \\
	
	Let $\lambda$ be any eigenvalue of $\rho(\sigma_1).$ Consider the corresponding eigenspace $E_{\lambda}\neq \{0\}.$ We  claim that it is invariant under $WB_n.$ Indeed, it is obviously invariant under $\rho(\sigma_i)=\rho(\sigma_1)$ for all
	$i=1, 2,\dots,n-1.$ And by the group identity (6), if $v\in E_{\lambda},$ then for any $i=1, 2,\dots,n-1,$ \\
	
	$\rho(\sigma_1)\left( \rho(\alpha_i)v\right)= \rho(\sigma_1)\left( \rho(\alpha_{n-1})v\right)=\rho(\alpha_{n-1}) \rho(\sigma_1)v=\lambda \left( \rho(\alpha_{n-1})v\right) =$\\
	
	$= \lambda \left( \rho(\alpha_{i})v\right) ,$ hence $\rho(\alpha_i)v\in E_{\lambda}.$\\
	
	If $E_{\lambda}=V,$ then $\rho(\sigma_i)=\lambda I$ for all $i=1, 2,\dots,n-1,$ and then $\rho$ is a direct sum of one-dimensional representations  $X_n(\lambda,k)$ of $WB_n,$  where $k \in \{-1,1\},$ hence $\rho$ is reducible.
	
	If $E_{\lambda}\neq V,$ then it is a proper invariant subspace of $V,$ hence $\rho$ is reducible as well.

\end{proof}

{\bf Theorem 6.2.} Suppose $n\geqslant 7.$   Let $\rho:WB_n\to GL_{n-1}(\mathbb{C})$ be an $(n-1)-$dimensional representation of the  welded braid group $WB_n,$  such that its restriction $\rho|_{B_n}$ onto the braid group $B_n$ is {\it irreducible}. Then there exists a unique quadruple $(y,k,z,\beta'),$ such that  $\rho$ is equivalent to $X_n(y,k)\otimes\beta',$ where $y,\ z \in\mathbb{C}^*,$ $P_n(z)=1+z+z^2+\dots +z^{n-1}\neq 0,$  $k\in \{-1,1\},$ and $\beta'$ is a representation of $WB_n,$\\  $\beta' \in \{\widetilde{\beta_n} (z),\widehat{\beta_n} (z)\},$ if $z\neq 1,$ and \\$\beta'=\widetilde{\beta_n} (1),$ if $z=1.$

\begin{proof} Since $\rho|_{B_n}$is irreducible $(n-1)-$dimensional representation,   then, by \cite{Formanek}, Theorem 23, there exists a unique pair $(y,z),$  $y,z\in\mathbb{C}^*,$  $P_n(z)\neq 0,$ such that $\rho|_{B_n}$ is equivalent to  $\chi_n(y)\otimes \beta_n(z),$ where $\beta_n(z)$ is the specialization of the reduced Burau representation, and $\chi_n(y)$ is a one-dimensional representation of $B_n.$ \\
	
Consider $\hat{\rho}=X_n(y^{-1},1)\otimes \rho.$ Its restriction $\hat{\rho}|_{B_n}$  onto the braid group $B_n$  is equivalent to $\beta_n(z),$ hence,  by Theorem 3.25, there exists a unique pair $(k,\beta'),$ such that  $\hat{\rho}$ is equivalent to $X_n(1,k)\otimes\beta',$ where $k\in \{ -1,1\},$ and $\beta'$ is a representation of $WB_n,$\\  $\beta'\in \{\widetilde{\beta_n} (z),\widehat{\beta_n} (z)\},$ if $z\neq 1,$ and \\$\beta'=\widetilde{\beta_n} (1),$ if $z=1.$\\ Thus,  $\rho=X_n(y,1)\otimes \hat{\rho}$ is equivalent to\\ $X_n(y,1)\otimes X_n(1,k)\otimes \beta'=X_n(y,k)\otimes  \beta',$ and the choice of $(y,k,z,\beta')$ is unique.\\
	
\end{proof}

{\bf Theorem 6.3.} Suppose $n\geqslant 9$  or $n=7.$ Let $\rho:WB_n\to GL_n(\mathbb{C})$ be an $n-$dimensional representation of the  welded braid group $WB_n,$  such that its restriction $\rho|_{B_n}$ onto the braid group $B_n$ is {\it irreducible}. Then there exists a unique quadruple $(y,k,z,\lambda),$ where $y,\ z,\ \lambda \in\mathbb{C}^*,$ $z\neq 1$ and $k\in \{-1,1\},$ such that $\rho$ is equivalent to $X_n(y,k)\otimes\widetilde{\tau}_n(z,\lambda).$

\begin{proof} The proof of this theorem is very similar to the proof of the Theorem 6.2. Since $\rho|_{B_n}$ is irreducible $n-$dimensional representation,  then for $n\geqslant 9$  by \cite{S}, Theorem 6.1, there exists a unique pair $(y,z),$  $y,z\in\mathbb{C}^*,$ $z\neq 1,$ such that  $\rho|_{B_n}$ is equivalent to  $\chi_n(y)\otimes \tau_n(z),$ where $\tau_n(z)$ is the specialization of the standard representation. The same statement for $n=7$ follows from  \cite{FLSV}, Theorem 6.1. \\

 Then the restriction of  $\hat{\rho}=X_n(y^{-1},1)\otimes \rho$   onto the braid group $B_n$  is equivalent to $\tau_n(z),$ $z\neq 1,$ and hence,  by Theorem 4.11, there exists a unique pair $(\lambda, k)$ where $\lambda \in\mathbb{C}^*$ and $k\in \{-1,1\}$ such that $\hat{\rho}$ is equivalent to $X_n(1,k)\otimes\widetilde{\tau}_n(z,\lambda).$ Hence $\rho=X_n(y,1)\otimes \hat{\rho}$ is equivalent to $X_n(y,1)\otimes X_n(1,k)\otimes \widetilde{\tau}_n(z,\lambda)=X_n(y,k)\otimes \widetilde{\tau}_n(z,\lambda),$ and the choice of  $(y,k,z,\lambda)$ is unique. \\

\end{proof}

\end{document}